\documentclass[AMA,Times1COL]{WileyNJDv5} 

\articletype{Research Article} 

\received{Date Month Year}
\revised{Date Month Year}
\accepted{Date Month Year}
\journal{Journal Name}
\volume{00}
\copyyear{2026}
\startpage{1}

\usepackage{booktabs}      
\usepackage{multirow}      
\usepackage{float}         
\usepackage{subcaption}    
\DeclareCaptionFont{intermediate}{\fontsize{9.5pt}{11pt}\selectfont}
\usepackage{algorithm}
\usepackage{algpseudocode}
\usepackage{amsmath}
\usepackage{bm}
\usepackage{hyperref}
\usepackage{xcolor}
\usepackage{amsfonts}
\usepackage{siunitx}

\begin{document}


\title{Physics informed learning of orthogonal features with applications in solving partial differential equations}

\author[1]{Qianxing Jia}

\author[1,2,3]{Dong Wang}

\authormark{JIA \textsc{et al.}} 
\titlemark{PHYSICS INFORMED LEARNING OF ORTHOGONAL FEATURES} 

\address[1]{\orgdiv{School of Science and Engineering}, \orgname{The Chinese University of Hong Kong (Shenzhen)}, \orgaddress{\state{Guangdong}, \postcode{518172}, \country{China}}}

\address[2]{\orgdiv{Shenzhen International Center for Industrial and Applied Mathematics}, \orgname{Shenzhen Research Institute of Big Data}, \orgaddress{\state{Guangdong}, \postcode{518172}, \country{China}}}

\address[3]{\orgname{Shenzhen Loop Area Institute}, \orgaddress{\state{Guangdong}, \postcode{518048}, \country{China}}}

\corres{Dong Wang, School of Science and Engineering, The Chinese University of Hong Kong (Shenzhen), Shenzhen, Guangdong 518172, China. \email{wangdong@cuhk.edu.cn}}

\abstract[Abstract]{The random feature method (RFM) constructs approximation spaces by initializing features from generic distributions, which provides universal approximation properties to solve general partial differential equations. However, such standard initializations lack awareness of the underlying physical laws and geometry, which limits approximation. In this work, we propose the Physics-Driven Orthogonal Feature Method (PD-OFM), a framework for constructing feature representations that are explicitly tailored to both the differential operator and the computational domain by pretraining features using physics-informed objectives together with orthogonality regularization. This pretraining strategy yields nearly orthogonal feature bases. We provide both theoretical and empirical evidence that physics-informed pretraining improves the approximation capability of the learned feature space. When employed to solve Helmholtz, Poisson, wave, and Navier–Stokes equations, the proposed method achieves residual errors 2–3 orders of magnitude lower than those of comparable methods. Furthermore, the orthogonality regularization improves transferability, enabling pretrained features to generalize effectively across different source terms and domain geometries for the same PDE.}

\keywords{random feature method, physics-informed neural network, orthogonality}

\jnlcitation{\cname{%
\author{Jia Q.}, and
\author{Wang D.}} (\cyear{2026}), 
\ctitle{Physics informed learning of orthogonal features with applications in solving partial differential equations}, \cjournal{Numerical Methods for Partial Differential Equations}, \cvol{2026;00:1--20}.}
\maketitle


\section{Introduction}

The numerical solution of partial differential equations is fundamental to a wide range of disciplines, including heat transfer, fluid dynamics, electromagnetism, and financial modeling. Classical numerical methods including finite element, finite difference, and spectral methods are effective for structured domains with well-defined boundary conditions. However, these methods are limited by the curse of dimensionality and require labor-intensive meshing for complex geometries. Physics-informed neural networks (PINN) \cite{Jiequn2018hipdedl,RAISSI2019686,SIRIGNANO20181339,SUKUMAR2022bdpinn, WANG2022whypinnsfailtotrain} has emerged as an alternative framework for solving PDEs using machine learning. PINN approximates the solution of a PDE using a neural network trained to minimize the residuals of the governing equations. This approach is mesh-free and suitable for high-dimensional problems and irregular domains. Despite these advantages, the ill-conditioned loss landscape of PINNs makes SGD converge slowly and often get stuck in local minima \cite{rathore2024pinns}. PINNs optimize the mean-squared residuals of PDEs evaluated at randomly sampled collocation points. Alternative approaches employ different loss formulations. Variational methods such as deep Ritz \cite{yu2018deep} and variational PINN \cite{KHARAZMI2021vpinn,URIARTE2025VPINNls} reformulate the PDEs via energy principles, often improving stability in high dimensions and reducing the order of differentiation thus reducing the cost of auto-differentiation \cite{baydin2018automatic}. Other methods adopt weak formulations, such as the weak adversarial network \cite{bao2024wancoweakadversarialnetworks,ZANG2020109409}, which has demonstrated robustness for problems with nonsmooth boundary conditions. In terms of architectural innovations, several studies explore modified network structures to improve accuracy. Kolomogorov-Arnold network introduces learnable univariate nonlinearities, allowing for better interpretability and approximation accuracy \cite{jacob2024spikansseparablephysicsinformedkolmogorovarnold,liu2025kankolmogorovarnoldnetworks,Toscano2025FromPINNs}. Tensor neural network design the network structure to be a combination of univariate functions for each dimension, which enables accurate quadrature-based integration \cite{WANG2025Tensorfp,Wang2024TensorNNpee,Yifan2024Tensor}. A distinct line of research focuses on decoupling the optimization of the network’s last layer from the rest of the parameters, thereby transforming the originally non-convex optimization problem into a convex one with respect to the last-layer coefficients. By interpreting the neural network as a linear combination of features from hidden layers, one can solve the final-layer coefficients via least squares. This idea has led to two main categories of methods: those using random initialization, such as the local extreme learning machine \cite{DONG2021elm,DONG2021pretrainingelm,DONG2022hyperelm}, the random feature method \cite{chen2024highprecisionrandomizediterativemethods, Jingrun2022rfm, liu2023random}, the randomized neural network \cite{deryck2025approximationtheoryapplicationsrandomized,Shang2024Randomized} and TransNet \cite{LU2025transnetdomain,ZHANG2024tarnsnet}, and those using trained neural network, including the subspace neural network \cite{lin2024adaptiveneuralnetworksubspace,song2025alternatelyoptimizedsnnmethodacoustic,xu2024subspacemethodbasedneural}, the Least-Squares VPINN \cite{URIARTE2025VPINNls}, and the Deepoly \cite{liu2025deepoly}. The advantage of training over random initialization remains unclear. This work provides both theoretical justification and empirical evidence that training the features, despite higher computational cost, leads to better accuracy. Inspired by classical spectral methods, which benefit from orthogonal basis functions with exponential convergence rates, this paper proposes a novel framework termed PD-OFM, short for physics-driven orthogonal feature method. Our method aims to learn operator-adaptive features that are orthogonal to each other. Penalizing the deviation of the features from orthogonal matrix during training guides the network to construct a nearly orthogonal feature space. After training, PDEs are solved via least-squares projection onto the space spanned by the learned features. Besides, this work introduces a projection error to quantify the distance between the spanned space of the learned features and the eigenspace of the differential operator. The projection error is reduced when orthogonality is regularized during training. Experiments on one-dimensional Poisson problems show that the features of PD-OFM closely approximates the analytical eigenfunctions and achieves high accuracy. Moreover, transferability helps reduce computational cost in repeated PDE solving. In artificial intelligence, it is common to freeze the pretrained feature space and then perform feature extraction or fine-tuning in subsequent tasks \cite{Bengio2021Deep, Pan2010Transfer}. This idea has later been successfully applied to solving partial differential equations \cite{LU2025transnetdomain, ZHANG2024tarnsnet, CHAKRABORTY2021Transfer, desai2022oneshot} and parametric partial differential equations \cite{CHEN2024104047,CHEN2024117198}. Furthermore, operator learning \cite{he2024mgno,li2021fourier,LiPINO2024,TRIPURA2023waveoperator} exhibit transferability, where the learned mappings are applied to new input functions or domains. The proposed PD-OFM also exhibits strong transferability. The pretrained features can also be transferred to solve the same type of PDE with different source terms or on varied domains. The main contributions of this paper are summarized as follows:
\begin{itemize}
    \item Theoretical analysis and numerical experiments demonstrate that physics-informed feature training enhances the approximation quality of the learned feature space, leading to improved accuracy in subsequent least-squares solution procedures.
    \item An orthogonality-regularized loss is introduced to ensure that the neural network’s penultimate-layer feature matrix yields a set of nearly orthogonal features.
    \item Comprehensive numerical experiments compare PD-OFM with random feature method and transferable neural network, highlighting its superior accuracy, stability, and transferability.
\end{itemize}
The remainder of this paper is organized as follows. Sections~\ref{subsec:review}–\ref{subsec:pretraining} review random feature method and develop the framework of physics-driven orthogonal feature method. Theoretical insights together with empirical evidence comparing trained and random representation functions are discussed in Section~\ref{subsec:training}. The PD-OFM is evaluated using metrics such as projection error and effective rank, as presented in Section \ref{subsec:orthogonality}. Section~\ref{subsec:transferability} explores the transferability of PD-OFM bases across different PDE instances. Numerical evaluations and comparisons with existing approaches are presented in Section~\ref{sec:experiments}.

\section{Physics-Driven Orthogonal Feature Method}
 
\subsection{Random Feature Method}
\label{subsec:review}
The random feature method for solving partial differential equations (PDEs) is built upon: (i) representing the solution of a PDE by a linear combination of random features, and (ii) enforcing the governing equations through the collocation method. The random features are typically constructed using randomly initialized shallow neural networks. Therefore, the method forms a natural integration of classical numerical techniques and modern machine learning approaches. Consider a general PDE of the form
\[
\begin{cases}
\mathcal{L} u(x) = f(x), & x \in \Omega,\\[6pt]
\mathcal{B} u(x) = g(x), & x \in \partial\Omega ,
\end{cases}
\]
where $x=(x_1,\cdots,x_d)^T$, and $\Omega$ is a bounded and connected domain. The random feature approximation assumes 
\[
u(x) = \sum_{i=1}^N c_i\, u_i(x),
\]
and determines the coefficients $\{c_i\}$ by minimizing a penalized least-squares loss:
\begin{equation}
    \mathcal{L}(c)
= \frac{1}{N_r} \sum_{i=1}^{N_r} \left| \mathcal{L} u(x_r^{(i)}) - f(x_r^{(i)}) \right|^2 + \frac{\lambda_b}{N_b} \sum_{j=1}^{N_b} \left| \mathcal{B} u(x_b^{(j)}) - g(x_b^{(j)}) \right|^2,
\label{eq:loss}
\end{equation}
where $\{x_r^{(i)}\}$ and $\{x_b^{(j)}\}$ are interior and boundary collocation points, respectively. Although this loss function has the same structure as that in physics-informed neural networks (PINNs), the optimization procedure differs significantly: all random features $u_i(x)$ are fixed, and only the coefficients $c_i$ are optimized. Consequently, gradient descent is unnecessary, and direct solvers such as linear or nonlinear least squares can be applied:
\[
c = \arg\min_c \mathcal{L}(c).
\]

Each random feature typically takes the form
\[
u_i(x) = \sigma(w_i x + b_i),
\]
where the parameters $(w_i,b_i)$ are sampled according to standard neural network initialization schemes such as uniform, Xavier, or Kaiming initialization \cite{he2023kaiminginitial, glorot2010understanding}. This corresponds to a shallow neural network in which only the output layer is trainable. Moreover, PDE solutions often exhibit strong local structures, and globally defined random features with standard neural initializations may lack sufficient approximation capability. To address this limitation, the random feature method introduces locally supported features through a partition of unity (PoU). To construct the PoU, a set of points $\{x_n\}_{n=1}^{M_p} \subset \Omega$ are sampled, serving as the center of a component in the partition. For each $n$, the normalized coordinate is constructed:
\begin{equation*}
    \tilde{x}_n=\frac{1}{r_n}(x-x_n), \quad n=1,\cdots,M_p,
\end{equation*}
where $r_n$ is preselected. When $d=1$, the PoU function is either of
\[
\psi_n^{a}(x) = \mathbb{I}_{-1 \le \tilde{x} < 1},
\]

\[
\psi_n^{b}(x) =
\begin{cases}
\dfrac{1 + \sin(2\pi \tilde{x})}{2}, & -\dfrac{5}{4} \le \tilde{x} < -\dfrac{3}{4}, \\[8pt]
1, & -\dfrac{3}{4} \le \tilde{x} < \dfrac{3}{4}, \\[8pt]
\dfrac{1 - \sin(2\pi \tilde{x})}{2}, & \dfrac{3}{4} \le \tilde{x} < \dfrac{5}{4}, \\[8pt]
0, & \text{otherwise}.
\end{cases}
\]
High-dimensional PoU is constructed using the tensor porduct of one-dimensional PoU functions $\psi_n(x)=\Pi_{k=1}^d \psi_n(x_k)$. Putting together, the solution of the PDE takes form:
\[
u_M(x) = \sum_{n=1}^{M_p} \psi_n(x) \sum_{j=1}^{J_n}c_{nj} u_{nj}(\tilde{x}_j).
\]
For the multiscale problem, the approximation incorporates both global and local components:
\[
u_M(x) = \sum_{i=1}^N c^g_i u^g_i(x) + \sum_{n=1}^{M_p} \psi_n(x) \sum_{j=1}^{J_n}c_{nj} u_{nj}(\tilde{x}_j).
\]

The random feature method is undoubtedly a powerful tool that has demonstrated remarkable success in solving a broad class of PDEs, including those with complex geometries, kinetic formulations, and multiscale structures. Despite its advantages, several important challenges remain.

First, although the random feature method reformulates PDEs into linear systems, the resulting matrices are frequently severely ill-conditioned, with condition numbers reaching as large as $10^{16}$. Such ill-conditioning significantly degrades the performance of iterative solvers. Even though the size of linear systems considered in the original studies were relatively small, reducing the condition number is essential for extending the method to larger-scale problems. Recent work by Chen and Tan~\cite{chen2024highprecisionrandomizediterativemethods} introduces randomized preconditioners to address this issue, while Shang et al.~\cite{Shang2024Randomized, SUN2024Randomized} adopts a Petrov-Galerkin and discontinuous Galerkin formulation rather than the standard collocation approach to improve stability. A second challenge lies in the initialization of random features. The ideas of utilizing neural network initialization techniques to initialize features were first brought about by Dong \cite{DONG2021elm, DONG2021pretrainingelm, DONG2022projection}. Thereafter, many literature focus on the initialization of features including Random feature method. Moreover, TransNet~\cite{ZHANG2024tarnsnet} introduces a more uniformly distributed initialization across the physical domain. Other approaches incorporate feature pretraining in order to be adaptive to different PDEs, including the subspace neural network~\cite{xu2024subspacemethodbasedneural}, the adaptive neural subspace method~\cite{lin2024adaptiveneuralnetworksubspace}, and Deepoly~\cite{liu2025deepoly}. Importantly, these two challenges are not mutually exclusive. In this work, we propose a framework that pretrains random features in a physics-informed manner, supplemented with orthogonality regularization, thereby simultaneously improving feature quality and alleviating ill-conditioning.

\subsection{Pretraining Features with Orthogonality Regularization}
\label{subsec:pretraining}
\begin{figure}[!h]
    \centering
    \includegraphics[width=1.0\textwidth]{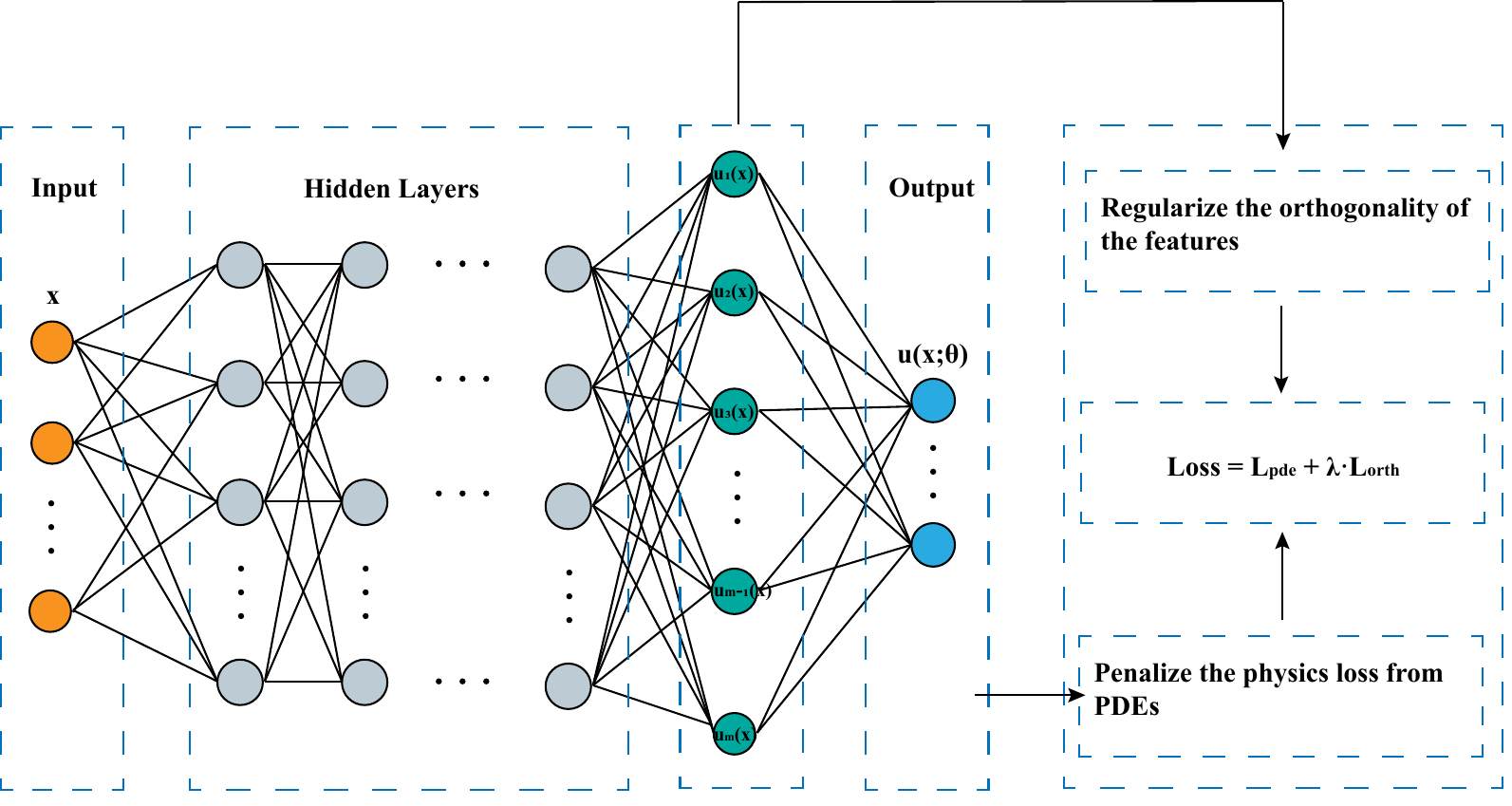}
    \caption{The network has two parallel output branches: one predicts the solution \(u(x; \theta)\), and the other outputs the feature space (green nodes).}
    \label{fig:network-architecture}
\end{figure}

A modified fully-connected neural network architecture is adopted, as illustrated in Figure~\ref{fig:network-architecture}. The network maps an input physical coordinate $x = (x_1, x_2, \dots, x_d) \in \Omega$ through a sequence of hidden layers. The network contains an intermediate feature space. The number of neurons, denoted by $m$ is treated as a hyperparameter controlling the number of features. The network produces two distinct outputs: the feature space layer outputs $$ U(x;\theta)= [ u_1(x;\theta),u_2(x;\theta), \cdots,u_m(x;\theta)]^T \in \mathbb{R}^m$$ and the final solution $u(x;\theta)=\sum_{i=1}^m c_i u_i(x;\theta)$, where $c_i$ is the coefficent of the last layer without bias. This output allows explicit regularization applied to the feature space layer.

The pretraining procedure is simultaneously minimizes all weights and biases by gradient descent with respect to the PDE loss and orthogonality regularizations:
\begin{equation}
\begin{aligned}
    \mathcal{L}_{\text{PINN}} &= \frac{1}{N_r} \sum_{i=1}^{N_r} \left| \mathcal{L} u(x_r^{(i)};\theta) - f(x_r^{(i)}) \right|^2 + \frac{1}{N_b} \sum_{j=1}^{N_b} \left| \mathcal{B} u(x_b^{(j)};\theta) - g(x_b^{(j)}) \right|^2\\
    \mathcal{L}_{\text{orth}} &= \left\| U^\top U - I \right\|_F^2\\
    \mathcal{L}_{\text{total}} &= \mathcal{L}_{\text{PINN}} + \lambda_{\text{orth}} \cdot \mathcal{L}_{\text{orth}}.
\end{aligned}
\end{equation}

\begin{algorithm}[!h]
\caption{PD-OFM: physics-driven orthogonal feature method}
\label{alg:PD-OFM}
\begin{algorithmic}[1]  
\Require Differential operator $\mathcal{L}$, boundary operator $\mathcal{B}$, source function $f$, boundary data $g$, collocation points $\{x_r^{(i)}\}$ and $\{x_b^{(j)}\}$, tolerance $\varepsilon$
\Ensure Approximate solution $u(x)$
\State Initialize neural network $u(x;\theta)$ with feature space layer of width $m$
\Repeat
    \State Sample interior and boundary collocation points
    \State Evaluate total loss: $\mathcal{L}_{\text{total}} = \mathcal{L}_{\text{PDE}} + \lambda_{\text{orth}} \cdot \mathcal{L}_{\text{orth}}$
    \State Update parameters $\phi$ via gradient descent
\Until{$\mathcal{L}_{\text{total}} < \varepsilon$ or convergence}
\State Freeze feature space layer output $U(x) = [u_1(x), \dots, u_m(x)]^T$
\State Form system matrix $A$ by evaluating $\mathcal{L} u_j(x_r^{(i)})$ and $\mathcal{B} u_j(x_b^{(j)})$ at all collocation points
\State Assemble right-hand side $b$ using $f(x_r^{(i)})$ and $g(x_b^{(j)})$
\State Solve the least-squares problem: $c^\ast = \arg\min_c \| A c - b \|_2^2$
\State \Return Approximate solution $u(x) = U(x)^\top c^\ast$
\end{algorithmic}
\end{algorithm}

\begin{remark}
When $\lambda_{\text{orth}} = 0$, Algorithm~\ref{alg:PD-OFM} reduces to the 
PD-FM (physics-driven feature method), which is equivalent to 
performing a re-optimization of the final linear coefficients after PINN training.
\end{remark}

\begin{figure}[t]
    \centering
    \begin{subfigure}[t]{0.45\linewidth}
        \centering
        \includegraphics[width=\linewidth]{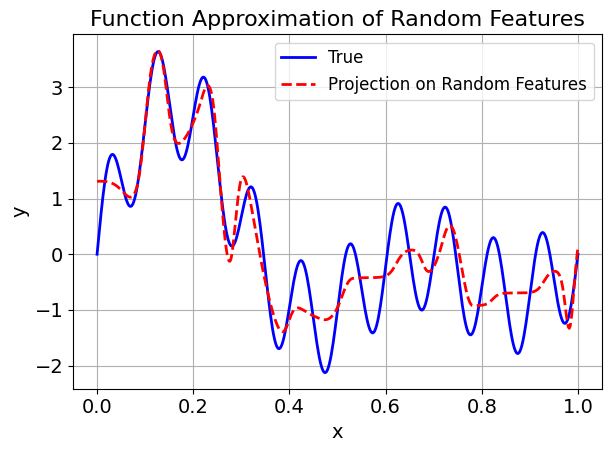}
        \caption{Projection of $f$ onto a random feature space.}
        \label{fig:approximate_rf}
    \end{subfigure}
    \hfill
    \begin{subfigure}[t]{0.45\linewidth}
        \centering
        \includegraphics[width=\linewidth]{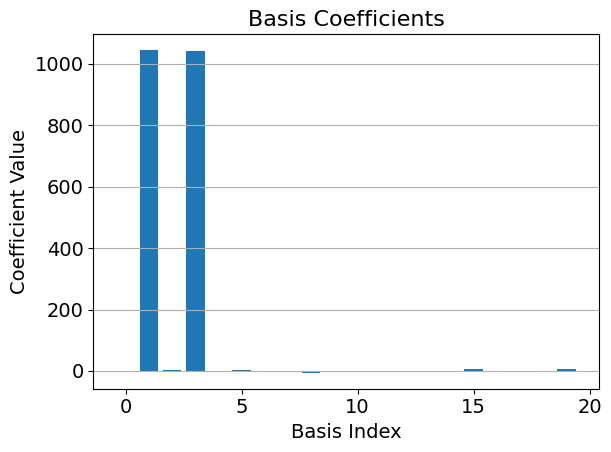}
        \caption{Distribution of coefficients associated with random features.}
        \label{fig:feature_coe}
   \end{subfigure}
    \caption{Approximation behavior of random features. The projection relies heavily on a small subset of features, indicating inefficient feature utilization.}
    \label{fig:rf_behavior}
\end{figure}

\begin{figure}[!h]
    \centering
    \includegraphics[width=0.45\linewidth]{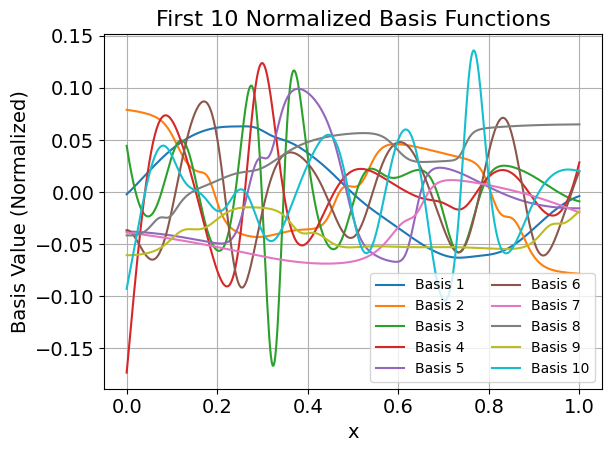}
    \caption{First 10 normalized basis functions learned by PD-OFM. The learned features remain well-distributed and non-degenerate, reflecting effective feature utilization.}
    \label{fig:feature_pdofm}
\end{figure}

\begin{figure}[!h]
    \centering
    \includegraphics[width=0.45\linewidth]{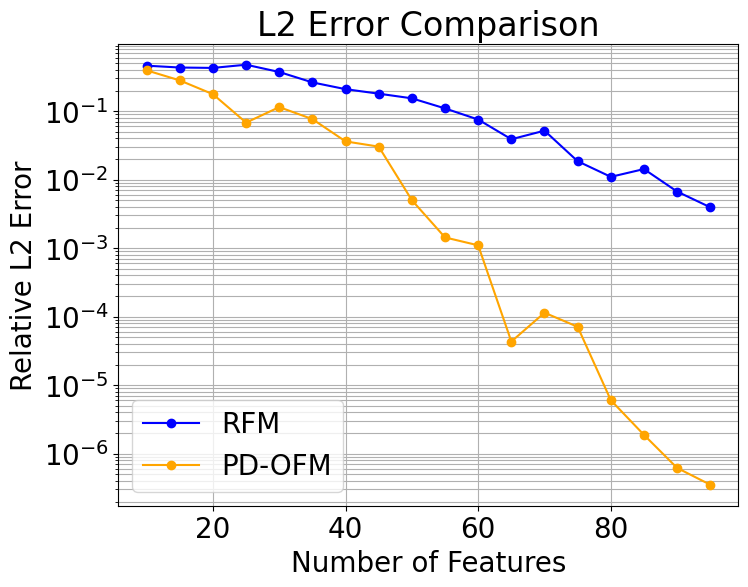}
    \caption{Approximation accuracy comparison between PD-OFM and the random feature method (RFM). PD-OFM exhibits consistently faster decay of the relative $L^2$ error as the number of features increases, achieving significantly higher accuracy at larger widths.}
    \label{fig:l2_error_pdofm}
\end{figure}

The motivation of PD-OFM stems from an intrinsic limitation of randomly initialized feature spaces in approximating functions with complex spectral content. In particular, when the number of features is small relative to the complexity of the target function, the approximation accuracy of random features deteriorates significantly. To illustrate this limitation, we consider the following one-dimensional function defined on $[0,1]$:
\[
f(x) = \sin(2\pi x) + \sin(3\pi x) + \sin(4\pi x) + \sin(20\pi x),
\]
which contains a mixture of low- and high-frequency components and therefore serves as a challenging test case for low-dimensional approximation spaces. To assess the intrinsic approximation capacity of random features, we project $f$ onto a randomly generated feature space with 20 features via a direct least-squares approximation in the function space. As shown in Fig.~\ref{fig:approximate_rf}, the resulting projection fails to capture the essential structure of the target function. Moreover, Fig.~\ref{fig:feature_coe} reveals that the optimal coefficients exhibit a highly uneven distribution, with a small subset of features dominating the approximation. In contrast, when the features are pretrained using a supervised loss together with an orthogonality regularization, a markedly different behavior emerges. As illustrated in Fig.~\ref{fig:feature_pdofm}, the learned features develop oscillatory patterns reminiscent of sinusoidal modes, suggesting that the training procedure induces a more expressive and structured feature space. Correspondingly, the approximation quality improves. As shown in Fig.~\ref{fig:l2_error_pdofm}, PD-OFM consistently achieves a faster decay of the relative $L^2$ error with respect to the number of features when compared to the random feature method.

This example provides evidence that physics-informed pretraining, coupled with orthogonality regularization, leads to more efficient use of finite-width representations.

\subsection{Training Does Matter}
\label{subsec:training}
In the preceding sections, the approximation behavior of different feature constructions was examined through static comparisons. The focus is now shifted to a dynamic perspective, with the objective of investigating how approximation ability changes throughout the training process. Earlier examples were formulated as supervised approximation problems, in which approximation quality was evaluated via direct least-squares fitting in the function space. In the context of solving partial differential equations, however, the approximation is governed by an operator-based loss, where the solution is obtained by minimizing the residual of a differential equation rather than matching a target function pointwise. This change in formulation motivates an investigation of training dynamics under operator-constrained settings. As a representative case, the one-dimensional Helmholtz equation is considered:
\begin{equation}
 \Delta u - k^2 u = f,
\end{equation}
with $k=\sqrt{10}$, where the source term $f$ is derived from the exact solution
\[
u(x) = \sin \left(3 \pi x + \frac{3 \pi}{20}\right)\cos\left(2 \pi x + \frac{\pi}{10}\right) + 2.
\]
Throughout this experiment, a fixed neural network architecture is employed, consisting of two hidden layers with 100 neurons. Residual connections \cite{He2016ResNet} are adopted together with the $\texttt{tanh}^3(\cdot)$ activation function \cite{bao2024wancoweakadversarialnetworks}. To examine the robustness of the observed training behavior, two widely used initialization strategies are considered: Xavier initialization \cite{glorot2010understanding} and Kaiming initialization \cite{he2023kaiminginitial}.

\begin{figure}[!h]
  \centering
  \begin{subfigure}[b]{0.42\textwidth}
    \centering
    \includegraphics[width=\textwidth]{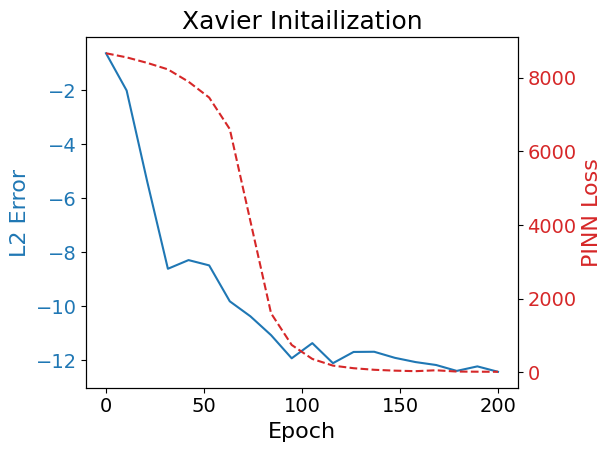}
    \label{fig:sub2}
  \end{subfigure}
  \hspace{0.02\textwidth}
  \begin{subfigure}[b]{0.4\textwidth}
    \centering
    \includegraphics[width=\textwidth]{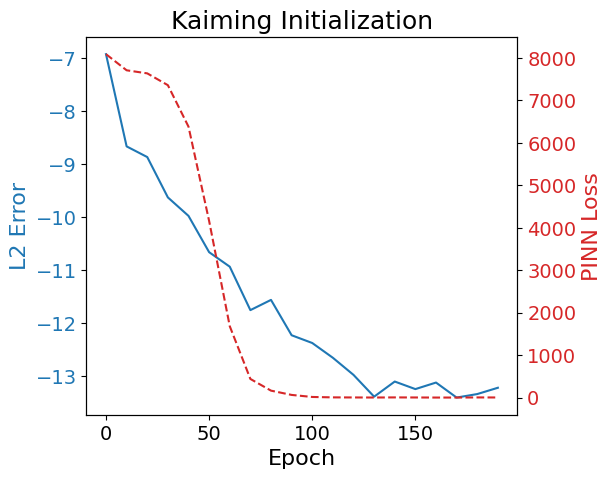}
    \label{fig:sub3}
  \end{subfigure}
  \caption{Training dynamics comparison between PINN loss and least-squares solution error.}
  \label{fig:doestraningmatter}
\end{figure}

As shown in Fig.~\ref{fig:doestraningmatter}, a clear correlation is observed between the evolution of the PINN residual loss and the accuracy of the least-squares solution. During the early stages of training, the decrease in the PINN loss is accompanied by a substantial reduction in the least-squares approximation error. As training progresses and the PINN loss approaches a plateau, further improvements in the least-squares solution become marginal, and the approximation error similarly stabilizes. Consistent behavior is observed across different initialization strategies. Without training, Xavier initialization leads to noticeably inferior approximation performance, whereas Kaiming initialization exhibits relatively better behavior. However, after only 100 training epochs, both initialization schemes experience significant improvements, and the discrepancy between them is largely reduced. This observation indicates that the training process progressively refines the feature space, enhancing the approximation capability of the learned features.

\textbf{Theorem 1} 
Let $A : L^2(\Omega) \to L^2(\Omega)$ be a bounded and stable linear operator, i.e., there exists a constant $\gamma > 0$ such that
\begin{equation}
    \|u\|_{L^2} \leq \gamma \|Au\|_{L^2}, \quad \forall u \in L^2(\Omega).
\end{equation}
Let $u_*$ be the exact solution to $Au = f$, and let $u_{NN}$ be the neural network approximation. Then,
\begin{equation}
    \|u_* - u_{NN}\|_{L^2} \leq \gamma \|Au_{NN} - f\|_{L^2}.
\end{equation}

\textit{Proof.}  
Since $A$ is stable, for any $v \in L^2(\Omega)$, we have $\|v\|_{L^2} \leq \gamma \|Av\|_{L^2}$. Letting $v = u_* - u_{NN}$, and using the linearity of $A$, we obtain
\begin{equation}
    \|u_* - u_{NN}\|_{L^2} \leq \gamma \|A(u_* - u_{NN})\|_{L^2} = \gamma \|f - Au_{NN}\|_{L^2}.
\end{equation}

When the PINN loss is sufficiently small, the bound indicates that the error between the neural network approximation \( u_{NN} \) and the exact solution \( u_* \) is also small. This suggests that the feature space learned by the network is capable of closely approximating the solution. Based on this feature space, a least-squares post-processing step can improve the solution accuracy.

\subsection{Why Orthogonality Matters: Structural Properties of the Feature Space}
\label{subsec:orthogonality}

While a low PINN training loss indicates that the network approximates the solution of a specific instance well, it does not guarantee the quality of the underlying feature space. To understand why PD-OFM outperforms standard methods, we must examine the structural properties of the learned basis functions. In classical spectral and Galerkin methods, the orthogonality of the basis is crucial for maximizing the approximation capacity of a subspace of fixed dimension and for ensuring the numerical stability of the linear solver. Motivated by this, we investigate three key metrics: \emph{effective rank}, \emph{orthogonality}, and \emph{projection error}. Our hypothesis is that orthogonality regularization forces the neural network to span a larger feature space, preventing feature collapse where multiple neurons learn redundant representations \cite{condensation2022Zhou}. To quantify these properties, we consider the benchmark one-dimensional Poisson equation:
\begin{equation}
    \begin{cases}
        -\dfrac{d^2 u(x)}{dx^2} = 1, & x \in (-1, 1), \\
        u(-1) = u(1) = 0.
    \end{cases}
\end{equation}
The analytical eigenfunctions of the Dirichlet Laplacian on $(-1, 1)$, $\phi_k(x) = \sin\left( \frac{k \pi (x + 1)}{2} \right)$, form a complete orthogonal basis. This allows us to directly measure the discrepancy between the learned feature space and the optimal spectral basis. We introduce the following definitions to quantify the quality of the feature space.

\begin{definition}
For a neural network \( u(x; \theta) \) with feature layer output \( \{u_j(x)\}_{j=1}^m \), the associated Gram matrix \( M_u \in \mathbb{R}^{m \times m} \) is defined as \( (M_u)_{ij} = \int_{\Omega} u_i(x) u_j(x) \, dx \). The \( \epsilon \)-effective rank of the feature space is defined as:
\[
r_\epsilon(M_u) := \left| \{ \lambda_k(M_u) : \lambda_k(M_u) > \epsilon \} \right|,
\]
where \( \lambda_k \) are the eigenvalues of \( M_u \).
\end{definition}

Yang et al. \cite{yang2025effectiverankstaircasephenomenon} proposed the effective rank as a robust measure of the linear independence of basis functions, filtering out numerical noise. A higher effective rank implies that the network utilizes its width more efficiently. To measure how well the learned space captures the true physics, we define the \textbf{Projection Error}.

\begin{definition} 
The projection error quantifies the ability of the learned features \( U(x) \) to represent the true eigenfunctions \( \{\phi_k(x)\} \). For the \( k \)-th eigenfunction, the relative projection error is:
\[
\epsilon_k^{\text{proj}} = \frac{ \| \mathcal{P}_U \phi_k - \phi_k \|_2 }{ \| \phi_k \|_2 }, \quad \text{where } \mathcal{P}_U \phi_k = U(U^\top U)^{-1} U^\top \phi_k \text{ is the orthogonal projection.}
\]
The average projection error is computed over the first \( M_1 \) dominant eigenfunctions: \( \epsilon^{\text{proj}} = \frac{1}{M_1} \sum_{k=1}^{M_1} \epsilon_k^{\text{proj}} \).
\end{definition}

\begin{figure}[!h]
    \centering
    \begin{subfigure}[b]{0.3\textwidth}
        \centering
        \includegraphics[width=\textwidth]{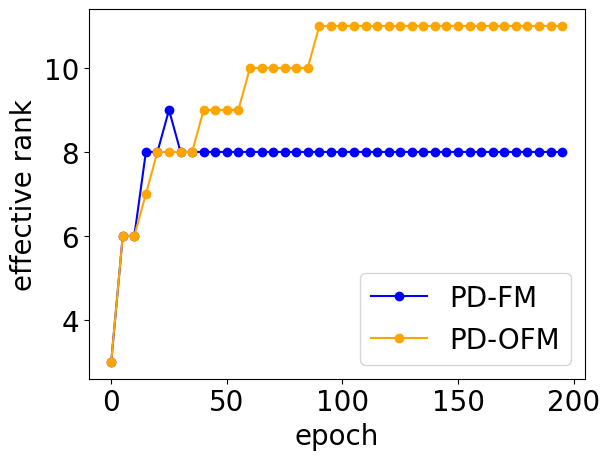}
        \caption{Effective rank}
        \label{fig:rank}
    \end{subfigure}
    \hfill
    \begin{subfigure}[b]{0.3\textwidth}
        \centering
        \includegraphics[width=\textwidth]{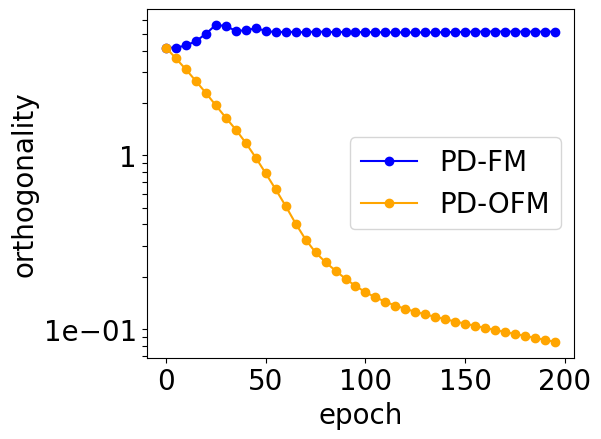}
        \caption{Orthogonality loss}
        \label{fig:ortho}
    \end{subfigure}
    \hfill
    \begin{subfigure}[b]{0.3\textwidth}
      \centering
        \includegraphics[width=\textwidth]{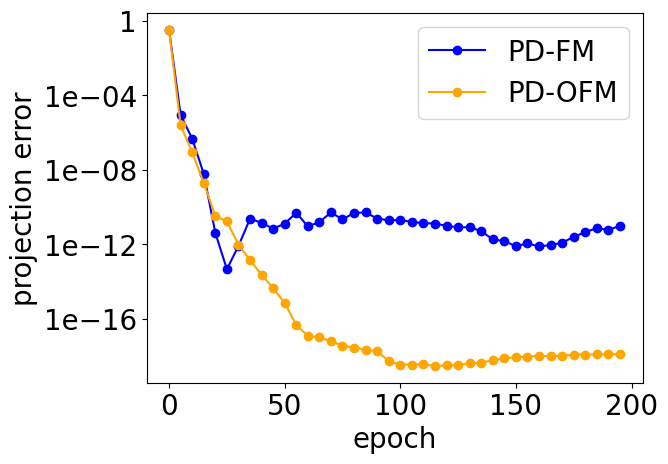}
        \caption{Projection error}
        \label{fig:proj}
    \end{subfigure}
    \caption{Evolution of feature space quality metrics during training. PD-OFM (orange) learns a higher-rank, more orthogonal, and more physically aligned basis compared to PD-FM (blue).}
    \label{fig:index_comp}
\end{figure}

This metric is rooted in the concept of Kolmogorov \(n\)-width \cite{pinkus2012n}, where the first \(k\) eigenfunctions span the optimal approximation subspace. We compare two models (with and without orthogonality regularization $\lambda_{\text{orth}} = 0.1$) trained under identical conditions using a ResNet architecture. Figure \ref{fig:index_comp} illustrates the evolution of these metrics:
\begin{itemize}
    \item \textbf{Effective Rank (\ref{fig:rank})}: The regularized model maintains a significantly higher rank, indicating that more neurons are contributing unique information.
    \item \textbf{Orthogonality (\ref{fig:ortho})}: The loss term effectively drives the feature covariance towards identity.
    \item \textbf{Projection Error (\ref{fig:proj})}: Crucially, the regularized features yield lower projection errors, meaning they spontaneously align better with the dominant eigenfunctions of the differential operator.
\end{itemize}

\begin{figure}[!h]
    \centering
    \includegraphics[width=0.5\linewidth]{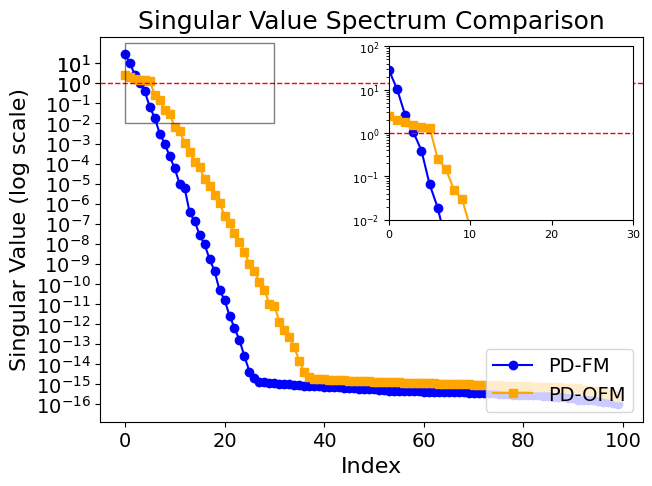}
    \caption{Singular value spectrum of the feature Gram matrix. PD-OFM exhibits a much slower decay rate than PD-FM.}
    \label{fig:sing}
\end{figure}

\subsubsection{Spectrum and Conditioning}
Figure \ref{fig:sing} presents the singular value spectrum of the feature Gram matrix. A key observation is that the singular values of PD-OFM decay much more slowly than those of PD-FM. While this does not yield a perfectly conditioned system (condition number $= 1$)—as the compactness of the inverse differential operator inevitably drives eigenvalues to zero—it significantly mitigates the severity of ill-conditioning. The spectrum remains flat for a larger number of indices before decaying. This slower decay implies that the system matrix in the least-squares step is numerically more stable.

\begin{figure}[!h]
    \centering
    \includegraphics[width=0.5\linewidth]{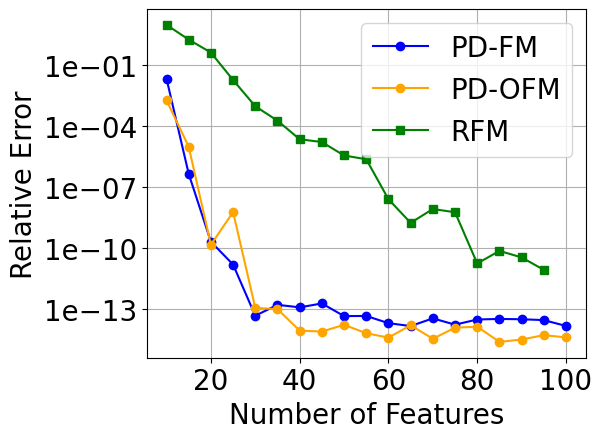}
    \caption{Relative $L^2$ error versus the number of features. The comparison reveals two advantages: (1) Training (PD-FM/OFM) outperforms random initialization (RFM); (2) Orthogonality (PD-OFM) further improves efficiency over standard training (PD-FM).}
    \label{fig:comp_feature}
\end{figure}

\subsubsection{Approximation Efficiency}
The impact of these structural properties on final accuracy is shown in Figure \ref{fig:comp_feature}, where we analyze the relative $L^2$ error as a function with respect to number of features. This result conveys two distinct layers of improvement:
\begin{enumerate}
    \item \textbf{Training vs. Random Initialization}: Both PD-FM and PD-OFM significantly outperform RFM. This confirms that physics-informed training successfully adapts the features to the specific operator, whereas random features may lack the necessary frequency components.
    \item \textbf{Orthogonality vs. Standard Training}: PD-OFM consistently achieves lower errors than PD-FM for the same number of features. This validates that orthogonality reduces redundancy; each added feature in PD-OFM contributes approximation power, whereas standard training suffers from diminishing returns due to feature correlation.
\end{enumerate}

\begin{remark}[Computational Complexity]
The overhead of orthogonality regularization is negligible. Let the network have $L$ layers of width $m$ and batch size $n$. Computing the regularization term involves matrix multiplications of cost $\mathcal{O}(n m^2)$. In contrast, the physics loss requires computing the Laplacian via reverse-mode automatic differentiation. For a problem in $d$ dimensions, this scales as $\mathcal{O}(d \cdot n L m^2)$ in the dominant term. Since typically $d, L \ge 2$, the gradient computation for the PDE residual dominates the cost, making the orthogonality constraint computationally free in relative terms.
\end{remark}

\subsection{Efficiency via Feature Transferability across Geometries}
\label{subsec:transferability}

The analysis in the previous section demonstrates that orthogonality regularization encourages the learned feature space to align with the spectral structure of the differential operator. Consequently, the pretrained features approximate the intrinsic eigenspace of the operator rather than overfitting to a specific solution instance. This property is pivotal for \textbf{computational efficiency}: it allows us to decouple the feature learning phase from the solution finding phase. Once a set of operator-adaptive features is learned, it can be reused as a universal basis to solve a wide class of problems sharing the same governing operator, significantly reducing the marginal cost of solving subsequent instances. To validate this efficiency, we propose a transfer learning framework where features are pretrained on a simple canonical domain and subsequently transferred to solve problems with varying source terms and, notably, different geometries. Consider the two-dimensional Poisson equation defined on a standard square domain \( \Omega_{base} = [-1, 1]^2 \):
\begin{equation}
-\frac{\partial^2 u}{\partial x^2} - \frac{\partial^2 u}{\partial y^2} = f, \quad (x, y) \in \Omega_{base}, \quad u|_{\partial \Omega_{base}} = 0.
\end{equation}
A feature set \( \{ u_i(x, y) \}_{i=1}^m \) is obtained by training PD-OFM on this base problem. These features are then frozen and transferred to solve new instances.

\subsubsection{Transfer across Source Terms}
We first evaluate the robustness of the features against variations in the source term $f(x,y)$. We construct a test dataset of 1000 randomized smooth solutions \( \{ u^{(j)} \}_{j=1}^{1000} \) generated by linear combinations of $K$ Gaussian kernels with random covariance structures:
\begin{equation}
u^{(j)}(x, y) = \sum_{k=1}^K c_k^{(j)} \exp \left( - \frac{1}{2} \mathbf{x}_k^\top (\Sigma_k^{(j)})^{-1} \mathbf{x}_k \right), \quad \mathbf{x}_k = \begin{bmatrix} x - \mu_{x,k}^{(j)} \\ y - \mu_{y,k}^{(j)} \end{bmatrix}.
\end{equation}
The coefficients are sampled as $c_k^{(j)} \sim \mathcal{N}(0, 1)$. The center vectors $\boldsymbol{\mu}_k^{(j)}$ and the covariance matrices $\Sigma_k^{(j)}$ are constructed as follows:
\begin{align}
\boldsymbol{\mu}_k^{(j)} &= \begin{bmatrix} \mu_{x,k}^{(j)} \\ \mu_{y,k}^{(j)} \end{bmatrix}, \quad \text{with} \quad \mu_{x,k}^{(j)}, \mu_{y,k}^{(j)} \sim \mathcal{U}[-1, 1], \\
\Sigma_k^{(j)} &= \begin{bmatrix} \sigma_{x,k}^2 & \rho_k \sigma_{x,k} \sigma_{y,k} \\ \rho_k \sigma_{x,k} \sigma_{y,k} & \sigma_{y,k}^2 \end{bmatrix}.
\end{align}
Here, the shape parameters are sampled uniformly as $\sigma_{x,k}, \sigma_{y,k} \sim \mathcal{U}[0.1, 0.5]$ and the correlation coefficient is $\rho_k \sim \mathcal{U}[-0.5, 0.5]$. The corresponding source terms $f_j$ are derived analytically.

\subsubsection{Transfer across Geometries}
A more challenging scenario involves transferring features to different geometric domains. Since different geometries impose different boundary conditions and imply different spectral bases, this test assesses whether the features learned on the superset domain \( \Omega_{base} \) typically capture the necessary frequency components to represent solutions on subsets. We transfer the features trained on the \textbf{Square} domain to an \textbf{Annulus} domain and an \textbf{L-shape} domain. For both scenarios, the approximate solution for the $j$-th instance is reconstructed via a computationally inexpensive least-squares fitting:
\begin{equation}
\hat{u}^{(j)}(x, y) = \sum_{i=1}^m c_{j,i} u_i(x, y),
\end{equation}
where only the linear coefficients \( c_{j,i} \) are re-computed to satisfy the PDE and boundary constraints of the new instance. The approximation quality is measured by the relative \( L^2 \) error:
\begin{equation}
\text{Err}_j = \frac{ \| \hat{u}^{(j)} - u^{(j)} \|_{L^2} }{ \| u^{(j)} \|_{L^2} }.
\end{equation}

Table \ref{tab:pdobnet_transfer} summarizes the average relative errors. PD-OFM achieves errors that are 1-2 orders of magnitude lower than the baseline PD-FM across all domains. Notably, the features generalize well even to the L-shape and Annulus domains without any fine-tuning of the neural network weights. This confirms that orthogonality regularization produces a robust, operator-adaptive basis that supports efficient solving of multi-geometry problems, offering a significant advantage over methods that require retraining for every geometric variation.

\begin{table}[!h]
    \centering
    \caption{Efficiency validation: Average relative $L^2$ errors under source term and domain transfer. Features trained on a square domain by PD-OFM generalize significantly better than PD-FM.}
    \begin{tabular}{ccccc}
        \toprule
        \textbf{Source Term} & \textbf{Model} & \textbf{Square} & \textbf{L-shape} & \textbf{Annular} \\
        \midrule
        \multirow{2}{*}{$2\pi^2\sin(\pi x)\sin(\pi y)$} & PD-OFM & $2.35 \times 10^{-8}$ & $1.13 \times 10^{-9}$ & $2.92 \times 10^{-9}$ \\
                         
                                & PD-FM  & $2.92 \times 10^{-7}$ & $1.06 \times 10^{-8}$ & $2.80 \times 10^{-7}$ \\
        \midrule
        \multirow{2}{*}{$(1-x^2)(1-y^2)$} & PD-OFM & $2.35 \times 10^{-8}$ & $8.32 \times 10^{-8}$ & $5.84 \times 10^{-10}$ \\
                      
                    & PD-FM  & $2.12 \times 10^{-7}$ & $2.15 \times 10^{-6}$ & $9.02 \times 10^{-8}$ \\
        \bottomrule
    \end{tabular}
    \label{tab:pdobnet_transfer}
\end{table}

Unlike one-shot transfer learning approaches that require joint training on multiple physics instances, PD-OFM constructs a transferable basis from a \emph{single} PDE instance. Detailed comparisons of different network depths and widths are provided in Appendix \ref{sec:comp} (Tables \ref{tab:exp1}--\ref{tab:exp3}).

\section{Numerical Experiments}
\label{sec:experiments}

A series of numerical experiments on different PDEs is conducted to further show the performance of PD-OFM, including linear, nonlinear, time-dependent, and multidimensional problems. The test cases also include domains ranging from simple intervals to complex geometries such as annuli and L-shaped regions. The goal is to evaluate accuracy and efficiency of PD-OFM, in comparison with several existing approaches, including the random feature method and the TransNet. Two standard error metrics: the relative $L^2$ error and the maximum pointwise $L^\infty$ error are computed over a set of uniformly sampled test points:
\begin{equation}
\begin{aligned}
\|e\|_{L^2} &= \sqrt{ \frac{\sum_{i=1}^N |u_\theta(x_i) - u^*(x_i)|^2}{\sum_{i=1}^N |u^*(x_i)|^2} }, \\
\|e\|_{L^\infty} &= \max_{1 \leq i \leq N} |u_\theta(x_i) - u^*(x_i)|.
\end{aligned}
\end{equation}

The effective rank and condition number of the system matrix encountered in the least-squares step are also recorded, which reflect the feature space's expressiveness and numerical stability.

\subsection{Benchmark on a Standard Problem: Helmholtz Equation}
The one-dimensional Helmholtz equation is first considered:
\begin{equation}
\begin{cases}
\dfrac{d^2 u}{dx^2} - \lambda u(x) = f(x), & x \in (0, 2), \\
u(0) = h_1, \quad u(2) = h_2.
\end{cases}
\end{equation}
$\lambda = 10$ and $f(x)$ is chosen such that the analytical solution is
\begin{equation}
    u^*(x) = \sin\left(3\pi x + \frac{3\pi}{20} \right) \cos\left(2\pi x + \frac{\pi}{10} \right) + 2.
\end{equation}

\begin{figure}[!h]
    \centering
    \begin{subfigure}{0.31\textwidth}
        \centering
        \includegraphics[width=\linewidth]{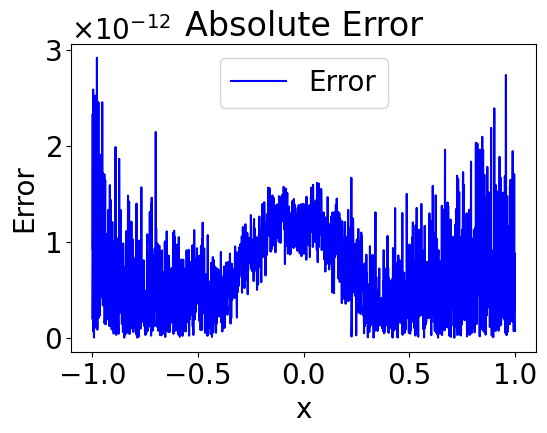}
        \caption{PD-OFM}
    \end{subfigure}
    \begin{subfigure}{0.3\textwidth}
        \centering
        \includegraphics[width=\linewidth]{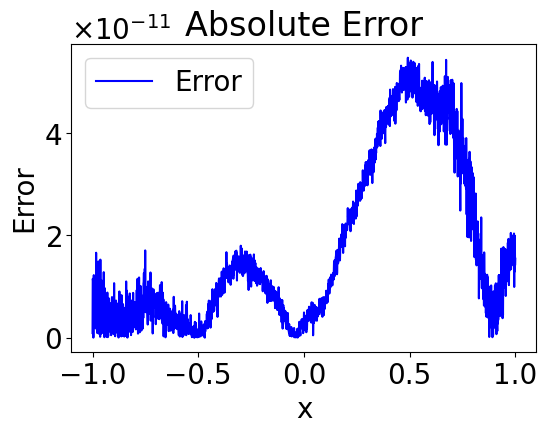}
        \caption{PD-FM}
    \end{subfigure}
    
    \begin{subfigure}{0.3\textwidth}
        \centering
        \includegraphics[width=\linewidth]{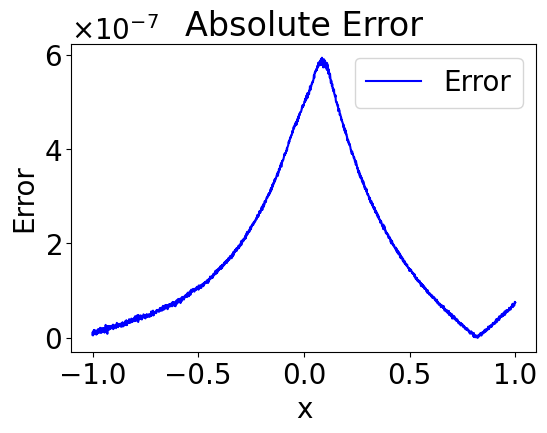}
        \caption{RFM}
    \end{subfigure}
    \begin{subfigure}{0.3\textwidth}
        \centering
        \includegraphics[width=\linewidth]{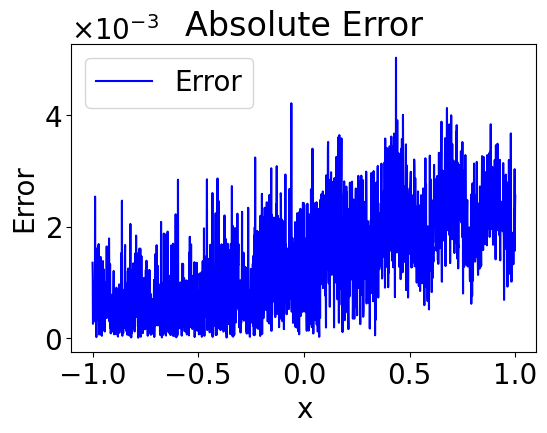}
        \caption{TransNet}
    \end{subfigure}

    \caption{Pointwise absolute errors for different methods on the 1D Helmholtz equation.}
    \label{fig:helmholtz_error}
\end{figure}

All models construct the solution as a linear combination of 100 features, corresponding to the width of the final hidden layer or the number of random features. For PD-OFM, the fully connected neural networks with \texttt{tanh}$^3$ activation functions is applied, trained using the Adam optimizer with initial learning rate $0.001$ and a maximum of $1000$ steps. The training is terminated early if the total loss falls below $\varepsilon = 10^{-3}$ of initial loss. For PD-OFM, the orthogonality regularization strength is set to $\lambda_{\text{orth}} = 1$. Interior collocation points are randomly sampled with a batch size of $1000$ at each iteration, and two fixed points are used to impose Dirichlet boundary conditions at the domain endpoints. The test set contains $2000$ equispaced sampled points over the domain. The results are shown in Table \ref{tab:helmholtz} and Figure \ref{fig:helmholtz_error}, where the proposed method reaches $10^{-13}$ accuracy.

\begin{table}[!h]
\centering
\begin{tabular}{lcccc}
\toprule
\textbf{Method} & \textbf{LS Residual} & \textbf{Rel. $L^2$ Error}   \\
\midrule
PD-OFM & $6.85 \times 10^{-10}$ & $3.98 \times 10^{-13}$   \\
PD-FM & $5.83 \times 10^{-7}$ & $1.38 \times 10^{-9}$   \\
RFM       & $1.01 \times 10^{-5}$  & $9.64 \times 10^{-9}$   \\
TransNet  & $1.00 \times 10^{-2}$  & $1.43 \times 10^{-3}$    \\
\bottomrule
\end{tabular}
\caption{Comparison of methods on the 1D Helmholtz equation. Errors are evaluated on 2000 uniformly spaced test points.}
\label{tab:helmholtz}
\end{table}

\subsection{Evaluating Complex Geometries: Poisson Equation on Irregular Domains}

The Poisson equation in two dimensions are considered:
\begin{equation}
\begin{cases}
-\Delta u(x,y) = f(x,y), & (x,y) \in \Omega, \\
u|_{\partial \Omega} = g(x,y).
\end{cases}
\end{equation}

To evaluate the robustness of each method with respect to domain geometry, three representative cases are considered:
\begin{itemize}
    \item \textbf{Square domain}: \(\Omega_1 = [-1, 1]^2\)
    \item \textbf{L-shaped domain}: \(\Omega_2 = [-1,1]^2 \setminus [0,1] \times [0,1]\)
    \item \textbf{Annular domain}: \(\Omega_3 = \{ (x,y) \in \mathbb{R}^2 \mid 0.25 \leq \sqrt{x^2 + y^2} \leq 1 \}\)
\end{itemize}

All models use 500 features. For PD-OFM, the neural network consists of two hidden layers of width 500 with \texttt{tanh}$^3$ activations, trained using the Adam optimizer with an initial learning rate of $0.001$ and maximum $1000$ steps. The training terminates when the total loss drops below $\varepsilon = 10^{-3}$. The orthogonality regularization strength is set to $\lambda_{\text{orth}} = 0.01$. In each case, $1024$ interior points and $128$ boundary points are sampled per iteration. Errors are evaluated on $2500$ uniformly spaced test points within the domain.

\begin{table}[!h]
\centering
\caption{Comparison of four methods on Poisson equation across different domains.}
\begin{tabular}{ccccc}
\toprule
Domain & Method & LS Residual & Rel.\ $L^2$ Error   \\
\midrule
Square Domain & PD-OFM & $2.87\times10^{-7}$ & $2.81\times10^{-9}$  \\
& PD-FM & $7.27 \times 10^{-7}$ & $9.62 \times 10^{-9}$  \\
         & RFM       & $7.83 \times 10^{-5}$ & $1.73 \times 10^{-5}$  \\
         & TransNet  & $1.47\times10^{-5}$  & $1.43\times10^{-7}$  \\
\midrule
L-shaped Domain & PD-OFM & $4.09 \times 10^{-8}$ & $6.48 \times 10^{-11}$  \\
& PD-FM 
& $5.98 \times 10^{-7}$ & $5.82 \times 10^{-9}$  \\
         & RFM & $2.88 \times 10^{-6}$ & $2.48 \times 10^{-6}$  \\
         & TransNet  & $1.04 \times 10^{-6}$ & $5.16 \times 10^{-9}$ \\
\midrule
Annulus Domain & PD-OFM & $5.43 \times 10^{-9}$ & $1.29 \times 10^{-11}$  \\
& PD-FM & $8.98 \times 10^{-9}$ & $4.13 \times 10^{-11}$  \\
         & RFM       & $4.93 \times 10^{-6}$ & $4.50 \times 10^{-8}$  \\
 
        & TransNet  & $4.30 \times 10^{-7}$ & $2.23 \times 10^{-9}$  \\
\bottomrule
\end{tabular}
\label{tab:poisson_3domains}
\end{table}

\begin{figure}[!h]
    \centering
    \begin{subfigure}{0.32\textwidth}
        \includegraphics[width=\linewidth]{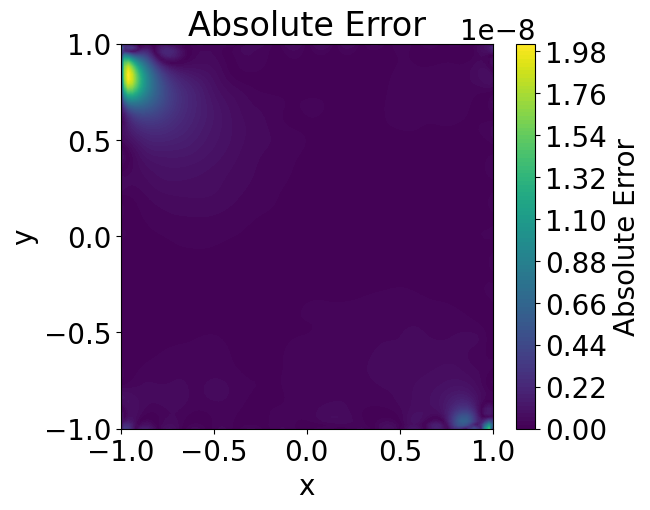}
        \caption{PD-OFM (Box)}
    \end{subfigure}
    \begin{subfigure}{0.32\textwidth}
        \includegraphics[width=\linewidth]{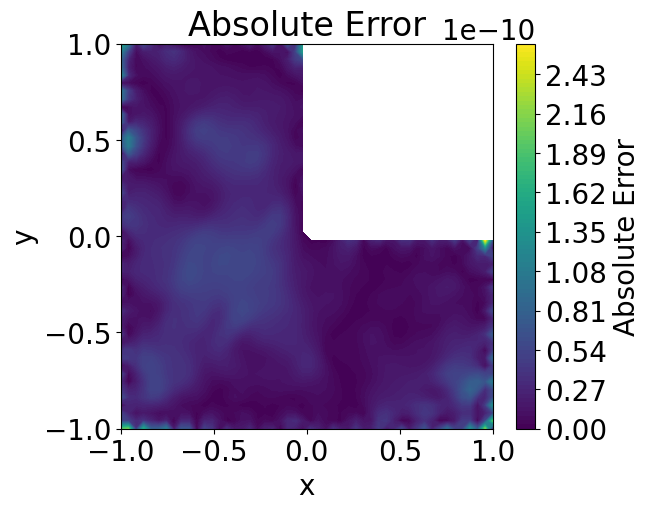}
        \caption{PD-OFM (L)}
    \end{subfigure}
    \begin{subfigure}{0.32\textwidth}
        \includegraphics[width=\linewidth]{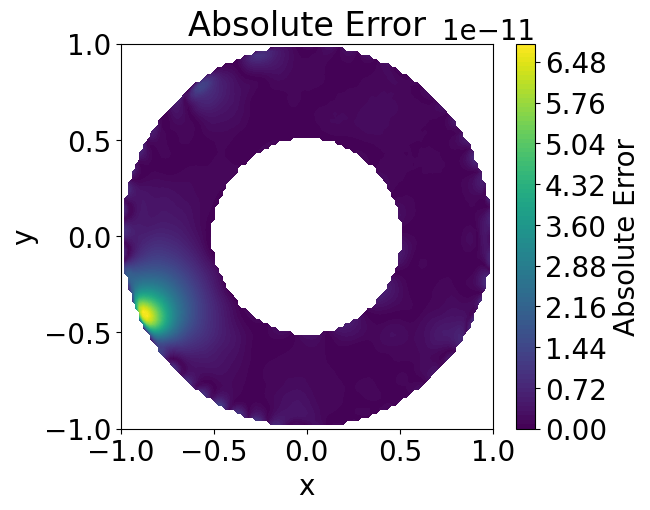}
       \caption{PD-OFM (A)}
    \end{subfigure}

    \begin{subfigure}{0.32\textwidth}
        \includegraphics[width=\linewidth]{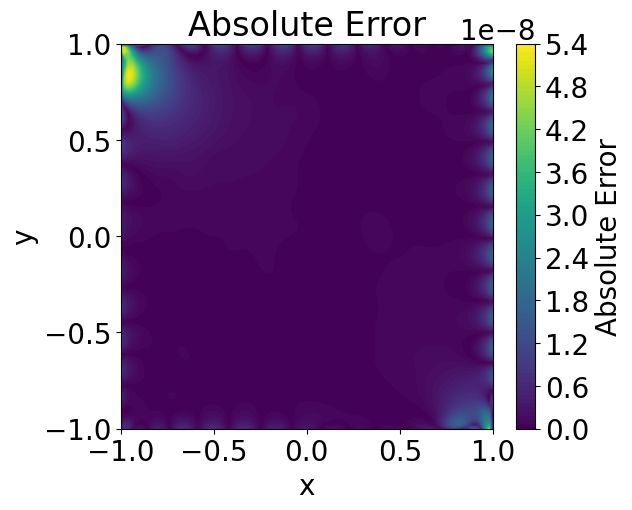}
        \caption{PD-FM (Box)}
    \end{subfigure}
    \begin{subfigure}{0.32\textwidth}
        \includegraphics[width=\linewidth]{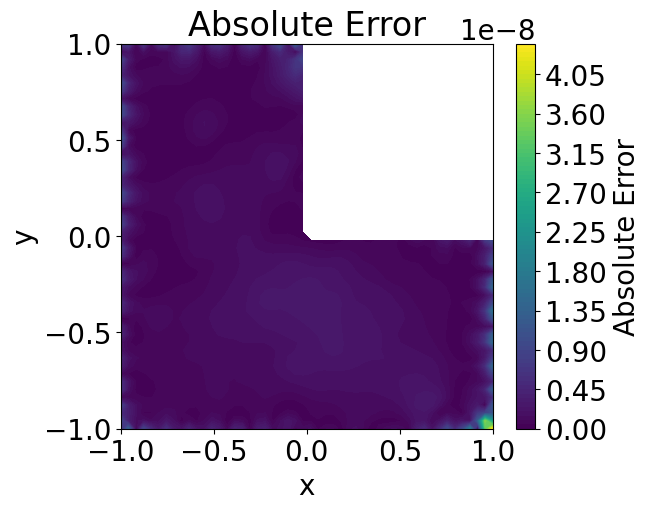}
        \caption{PD-FM (L)}
    \end{subfigure}
    \begin{subfigure}{0.32\textwidth}
        \includegraphics[width=\linewidth]{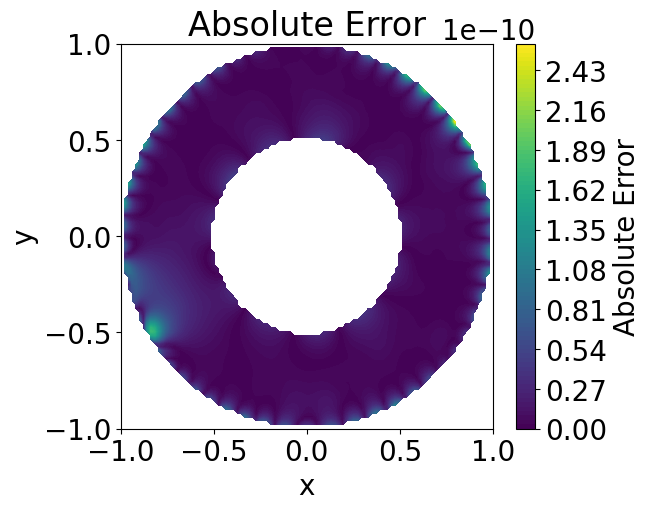}
        \caption{PD-FM (A)}
    \end{subfigure}

    \begin{subfigure}{0.32\textwidth}
        \includegraphics[width=\linewidth]{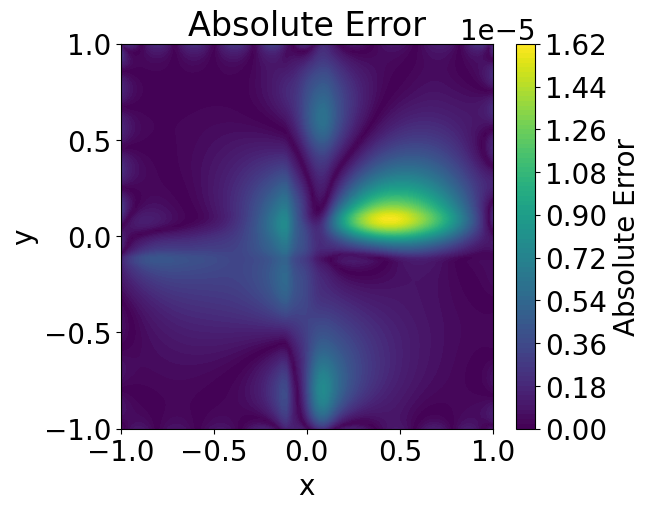}
        \caption{RFM (Box)}
    \end{subfigure}
    \begin{subfigure}{0.32\textwidth}
        \includegraphics[width=\linewidth]{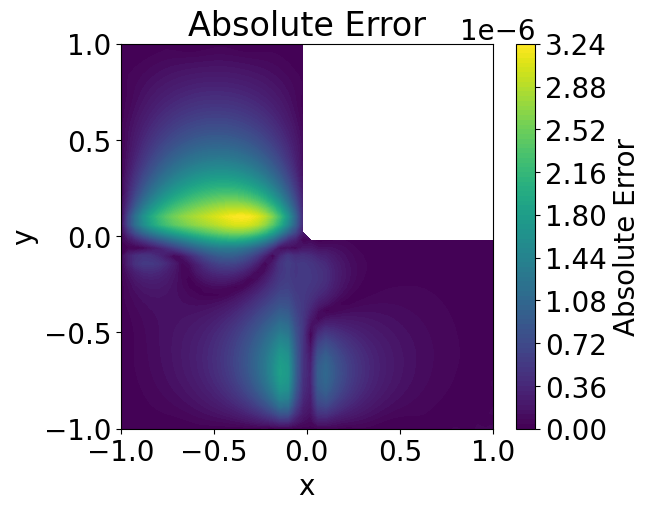}
        \caption{RFM (L)}
    \end{subfigure}
    \begin{subfigure}{0.32\textwidth}
        \includegraphics[width=\linewidth]{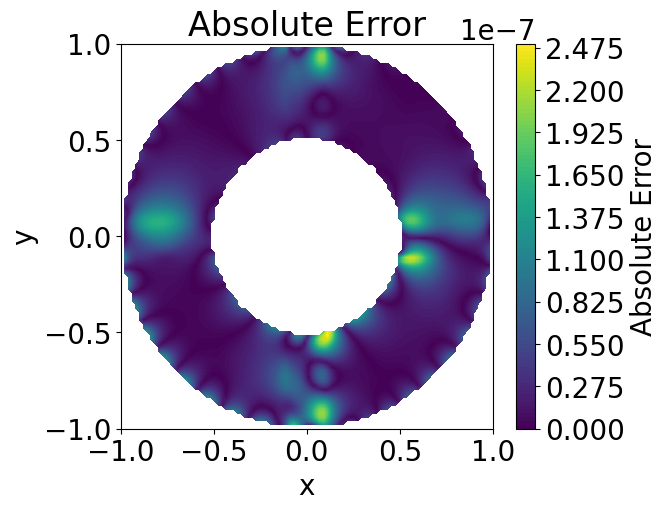}
        \caption{RFM (A)}
    \end{subfigure}

    \begin{subfigure}{0.32\textwidth}
        \includegraphics[width=\linewidth]{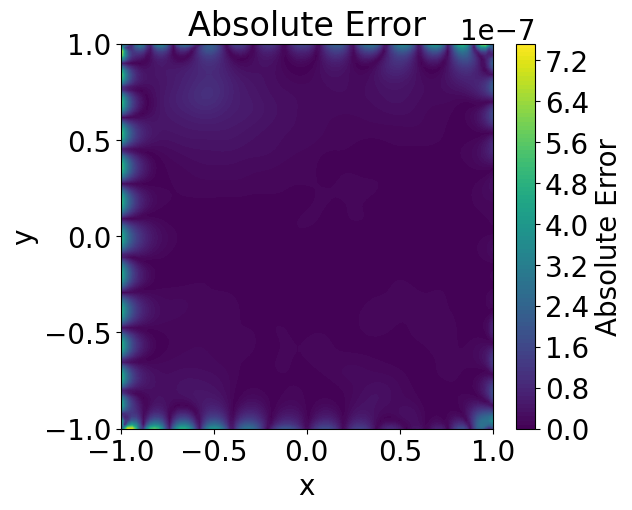}
        \caption{TransNet (Box)}
     \end{subfigure}
    \begin{subfigure}{0.32\textwidth}
        \includegraphics[width=\linewidth]{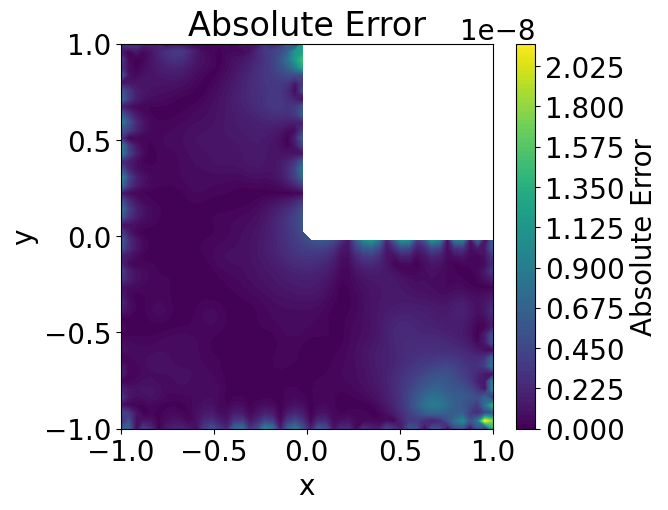}
        \caption{TransNet (L)}
    \end{subfigure}
    \begin{subfigure}{0.32\textwidth}
        \includegraphics[width=\linewidth]{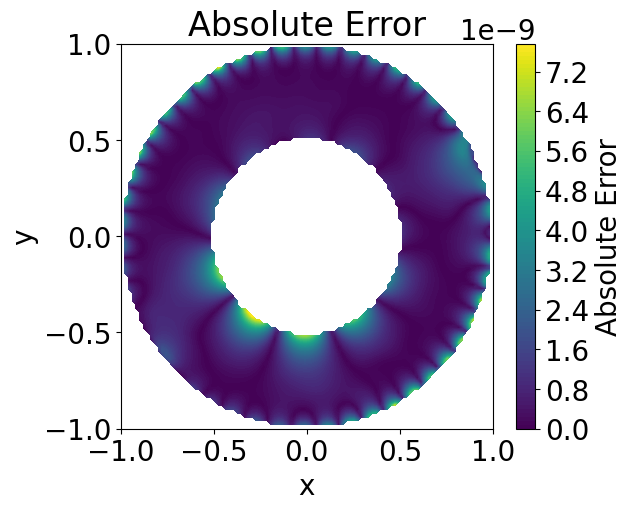}
        \caption{TransNet (A)}
    \end{subfigure}

    \caption{
    Pointwise absolute errors of different methods (rows) on the two-dimensional Poisson equation across three domains (columns): Box, L-shape, and Annulus.
}
    \label{fig:poisson_2d_domains}
\end{figure}

As shown in Table \ref{tab:poisson_3domains} and Figure \ref{fig:poisson_2d_domains}, PD-OFM consistently achieves the highest accuracy across all tested domains. Interestingly, each method exhibits distinct error patterns in their final predictions. In PD-OFM, the error is highly localized and accumulates at a single spatial point. This may be attributed to the orthogonality constraint encouraging globally balanced representations, while subtle local features are neglected. In contrast, the random feature method, despite using smooth partition of unity, shows the largest error at the interfaces between subdomains. This could be due to the inherent asymmetry between interior and boundary constraints in the least-squares formulation, which tends to under-emphasize boundary enforcement.

\subsection{Scalability to Higher Dimensions: three-dimensional Poisson Equation on Unit Cube}
Consider the three-dimensional Poisson equation on the domain $\Omega = [0, 1]^3$:
\begin{equation}
\begin{cases}
-\Delta u(x, y, z) = f(x, y, z), & (x, y, z) \in \Omega, \\
u(x, y, z) = 0, & (x, y, z) \in \partial \Omega.
\end{cases}
\end{equation}

The exact solution is given by:
\begin{equation}
u(x, y, z) = \sin \left( \pi x \right) \sin \left(\pi y \right) \sin \left( \pi z  \right),
\end{equation}
and the corresponding source term is computed analytically. The network architecture consists of two hidden layers with 600 neurons per layer and residual connections. The feature space layer width is set to 600. The orthogonality regularization strength is $\lambda_{\text{orth}} = 0.01$. A total of 2048 interior collocation points and 6,00 boundary collocation points are used.

\begin{table}[!h]
\centering
\caption{Comparison of methods on the three-dimensional Poisson equation. Errors are evaluted on 2700 uniformly spaced test points and taken average over z axis.}
\begin{tabular}{lcccc}
\toprule
\textbf{Method} & \textbf{LS Residual} & \textbf{Rel. $L^2$ Error}  \\
\midrule
PD-OFM & $6.89 \times 10^{-5}$ & $1.64 \times 10^{-5}$ \\
PD-FM & $6.37 \times 10^{-5}$ & $1.07 \times 10^{-4}$ \\
RFM       & $1.01 \times 10^{-2}$  & $1.38 \times 10^{-1}$  \\
TransNet  & $1.42 \times 10^{-4}$  & $2.03 \times 10^{-4}$   \\
\bottomrule
\end{tabular}
\label{tab:3dpoisson}
\end{table}

\begin{figure}[!h]
    \centering
    \begin{subfigure}{0.4\textwidth}
        \centering
        \includegraphics[width=\linewidth]{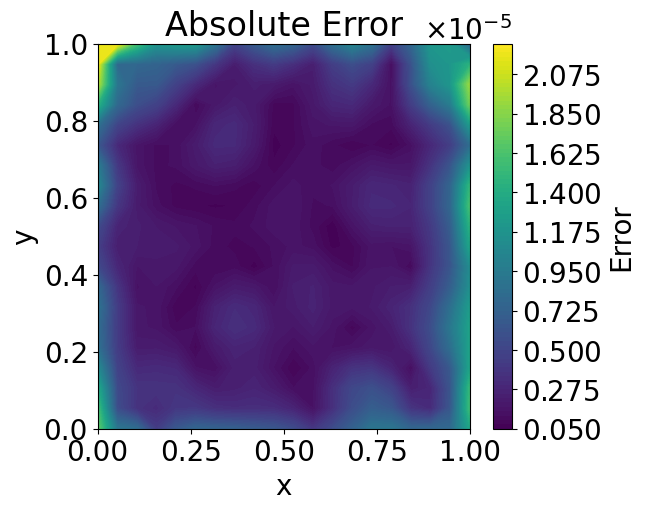}
        \caption{PD-OFM}
    \end{subfigure}
    \hspace{0.1\textwidth}
    \begin{subfigure}{0.4\textwidth}
        \centering
        \includegraphics[width=\linewidth]{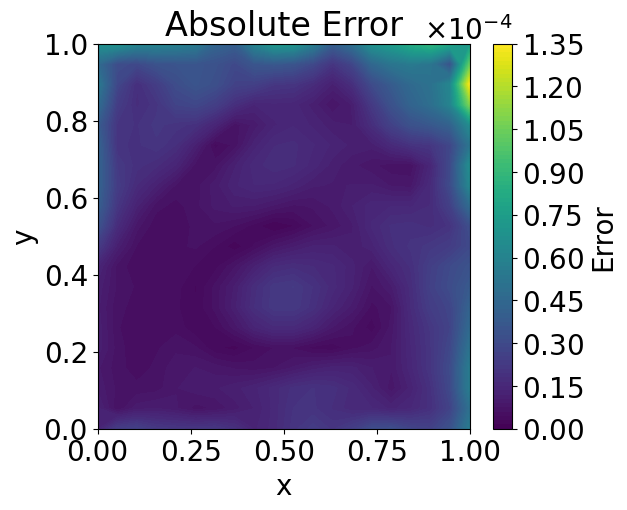}
        \caption{PD-FM}
    \end{subfigure}
    
    \begin{subfigure}{0.4\textwidth}
        \centering
        \includegraphics[width=\linewidth]{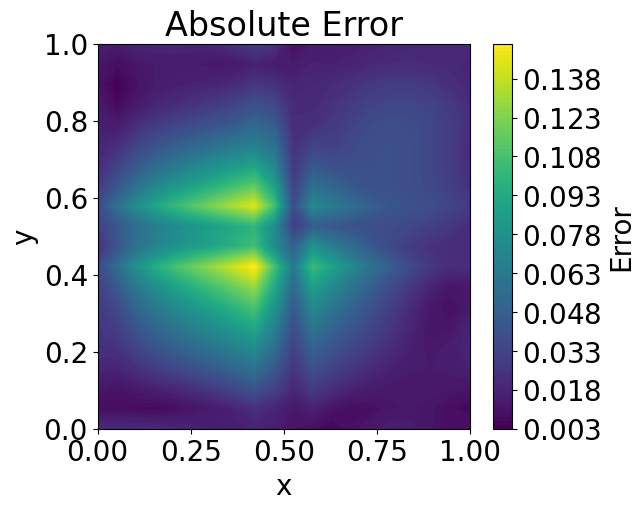}
        \caption{RFM}
    \end{subfigure}
    \hspace{0.1\textwidth}
    \begin{subfigure}{0.4\textwidth}
        \centering
        \includegraphics[width=\linewidth]{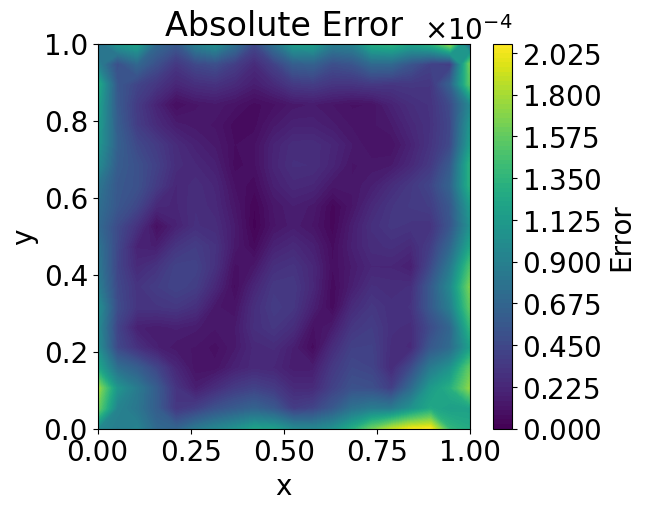}
        \caption{TransNet}
    \end{subfigure}

    \caption{Pointwise absolute errors for different methods on the three-dimensional Poisson equation.}
    \label{fig:3dpoisson_error}
\end{figure}

Table \ref{tab:3dpoisson} reports the quantitative comparison on the three-dimensional Poisson equation. PD-OFM achieves the lowest least-squares residual and relative $L^2$ error among all methods, with a moderate effective rank. Although RFM attains the highest rank due to the partition of unity, the learned features fails to capture the spectral structure of the PDE, as indicated by its significantly higher residual and projection error. The visualized pointwise absolute errors in Figure \ref{fig:3dpoisson_error} confirm that the RFM features, despite its large rank, cannot approximate the solution accurately. The TransNet method shows improved error compared to RFM but remains inferior to PD-OFM.

\subsection{Exploring Diverse Boundary Conditions: Wave Equation}

In this third test case, the time-dependent one-dimensional wave equation is considered:
\begin{equation}
\begin{cases}
\displaystyle \frac{\partial^2 u}{\partial t^2} = c \frac{\partial^2 u}{\partial x^2}, & x \in [0, 1],\ t \in [0, 2], \\[1.2ex]
u(x, 0) = \sin(4\pi x), & \\[1.2ex]
\displaystyle \left. \frac{\partial u}{\partial t} \right|_{t=0} = 0, & \\[1.2ex]
u(0, t) = u(1, t), &
\end{cases}
\end{equation}
where \( c = \frac{1}{16\pi^2} \). The exact solution to this equation is given by
\begin{equation}
    u(x,t)=\frac{1}{2}(\sin(4\pi x + t) - \sin(4\pi x - t)).
\end{equation}

\begin{figure}[!h]
    \centering
    \begin{subfigure}{0.4\textwidth}
        \centering
        \includegraphics[width=\linewidth]{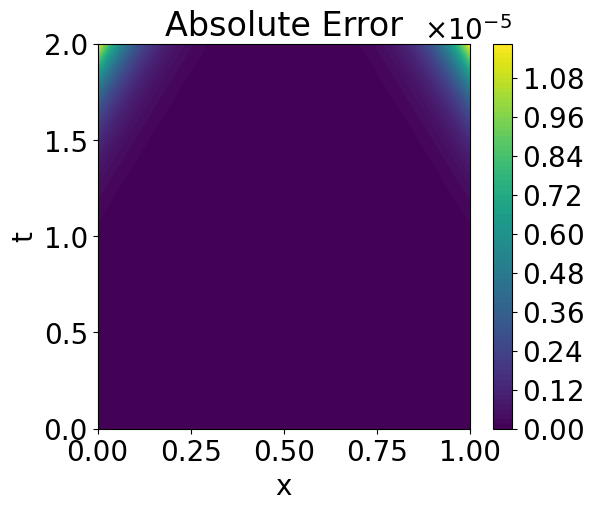}
        \caption{PD-OFM}
    \end{subfigure}
    \hspace{0.1 \textwidth}
    \begin{subfigure}{0.4\textwidth}
        \centering
        \includegraphics[width=\linewidth]{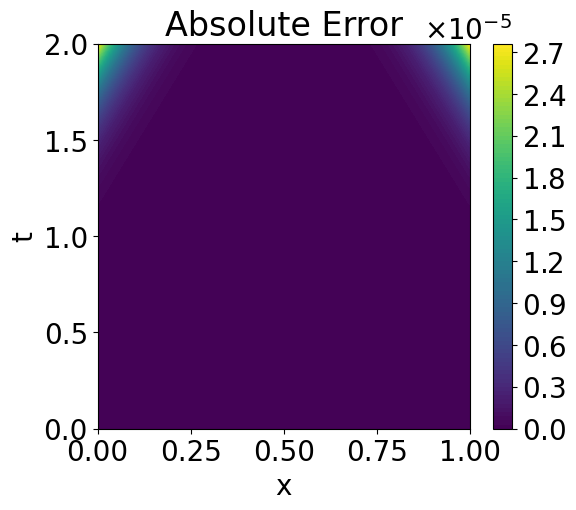}
        \caption{PD-FM}
    \end{subfigure}
    
    \begin{subfigure}{0.4\textwidth}
        \centering
         \includegraphics[width=\linewidth]{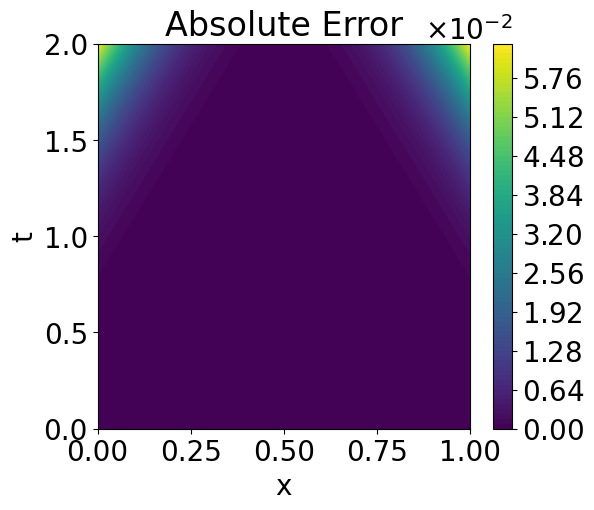}
        \caption{RFM}
    \end{subfigure}
    \hspace{0.1 \textwidth}
    \begin{subfigure}{0.4\textwidth}
        \centering
        \includegraphics[width=\linewidth]{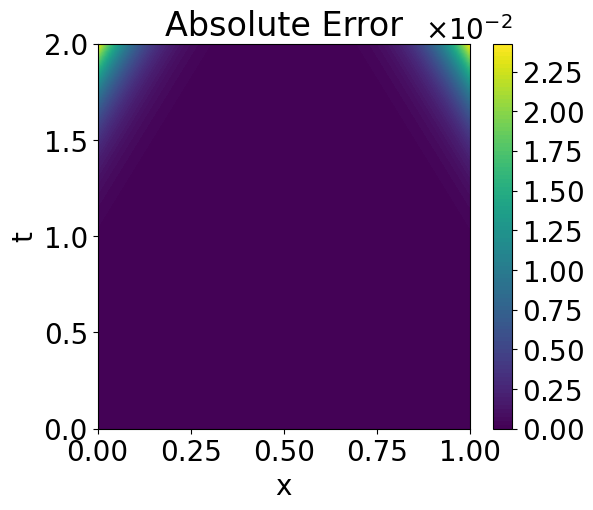}
        \caption{TransNet}
    \end{subfigure}
    \caption{Pointwise absolute errors for different methods on the 1D wave equation.}
    \label{fig:wave_error}
\end{figure}

This equation features not only the usual initial conditions (initial displacement and initial velocity), but also a periodic boundary condition. Our method is flexible and can be adapted to different types of boundary conditions. Specifically, the periodic boundary condition is enforced at points 
\[
\{(0,t_1),\dots,(0,t_n),\ (1,t_1),\dots,(1,t_n)\},
\]
which lie on the boundaries \( x=0 \) and \( x=1 \) and share the same set of time values \( \{t_1,\dots,t_n\} \). The results are shown in Table \ref{tab:wave_comp} and Figrue \ref{fig:wave_error}.
\begin{table}[!h]
\centering
\begin{tabular}{lccc}
\toprule
\textbf{Method} & \textbf{LS Residual} & \textbf{Rel. $L^2$ Error} \\
\midrule
PD-OFM & $6.45 \times 10^{-11}$ & $1.99 \times 10^{-6}$  \\
PD-FM & $1.25 \times 10^{-11}$ & $4.77 \times 10^{-6}$   \\
RFM       & $1.09 \times 10^{-6}$  & $1.62 \times 10^{-2}$    \\
TransNet  & $2.75 \times 10^{-9}$  & $4.41 \times 10^{-3}$  \\
\bottomrule
\end{tabular}
\caption{Comparison of methods on the 1D Wave equation. Errors are evaluated on 2500 uniformly spaced test points.}
\label{tab:wave_comp}
\end{table}

As shown in Figure~\ref{fig:wave_error}, all methods yield larger errors compared to the previous examples. Notably, when the periodic boundary condition is replaced with a Dirichlet boundary condition, the absolute error drops below \( 1 \times 10^{-10} \). This suggests that the periodic boundary condition introduces an imbalance when combined with the other constraints. Nevertheless, PD-OFM achieves the lowest error among all the tested methods.

\subsection{Nonlinear PDEs: Steady-State Navier–Stokes}

The steady-state incompressible Navier--Stokes equations are considered:
\begin{equation}
\begin{cases}
\boldsymbol{u} \cdot \nabla \boldsymbol{u} + \nabla p - \nu \Delta \boldsymbol{u} = 0, \\
\nabla \cdot \boldsymbol{u} = 0,
\end{cases}
\end{equation}
where $\bm{u} = (v_1, v_2)$ is the velocity field, $p$ is the pressure, and $\nu$ denotes the kinematic viscosity. The Reynolds number is defined as $Re = 1/\nu$.

As a benchmark, the classical Kovasznay flow is adopted, which provides an exact analytical solution to the steady Navier--Stokes equations in two dimensions:
\begin{align}
v_1(x_1, x_2) &= 1 - e^{\lambda x_1} \cos(2\pi x_2), \\
v_2(x_1, x_2) &= \frac{\lambda}{2\pi} e^{\lambda x_1} \sin(2\pi x_2), \\
p(x_1, x_2) &= \frac{1}{2} \left(1 - e^{2\lambda x_1} \right),
\end{align}
where
\begin{equation}
    \lambda = \frac{1}{2\nu} - \sqrt{\frac{1}{4\nu^2} + 4\pi^2}.
\end{equation}

\begin{figure}[!h]
    \centering

    \begin{subfigure}[t]{0.3\linewidth}
        \centering
        \includegraphics[width=\linewidth]{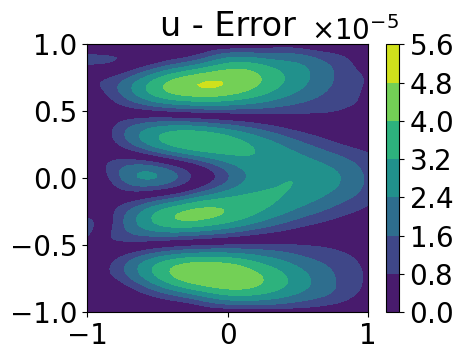}
        \caption{PD-OFM (u)}
    \end{subfigure}
    \hfill
    \begin{subfigure}[t]{0.3\linewidth}
        \centering
        \includegraphics[width=\linewidth]{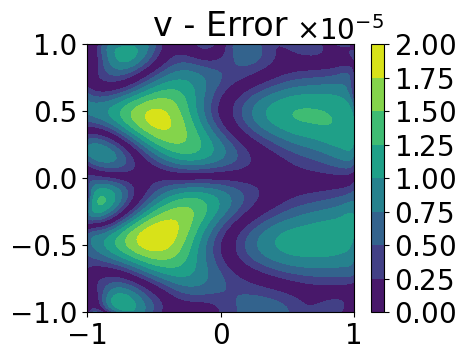}
        \caption{PD-OFM (v)}
    \end{subfigure}
    \hfill
    \begin{subfigure}[t]{0.3\linewidth}
        \centering
        \includegraphics[width=\linewidth]{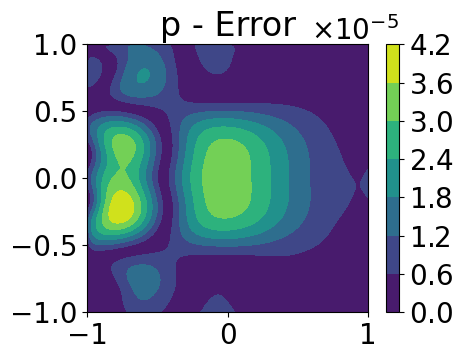}
        \caption{PD-OFM (p)}
    \end{subfigure}

    \begin{subfigure}[t]{0.3\linewidth}
        \centering
        \includegraphics[width=\linewidth]{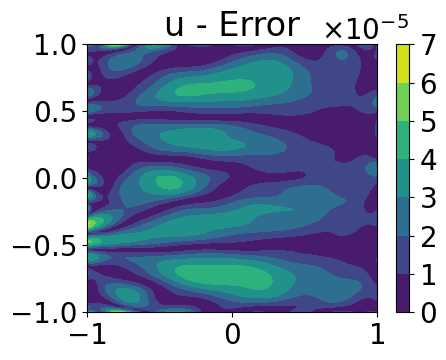}
        \caption{PD-FM (u)}
    \end{subfigure}
    \hfill
    \begin{subfigure}[t]{0.3\linewidth}
        \centering
        \includegraphics[width=\linewidth]{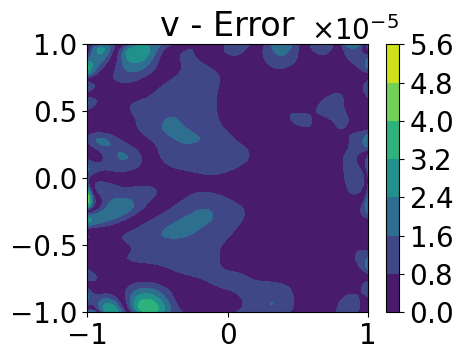}
        \caption{PD-FM (v)}
    \end{subfigure}
    \hfill
    \begin{subfigure}[t]{0.3\linewidth}
        \centering
        \includegraphics[width=\linewidth]{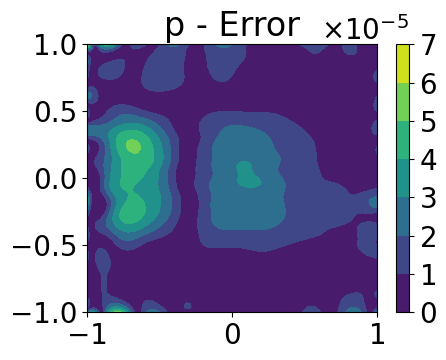}
        \caption{PD-FM (p)}
    \end{subfigure}

    \begin{subfigure}[t]{0.3\linewidth}
        \centering
        \includegraphics[width=\linewidth]{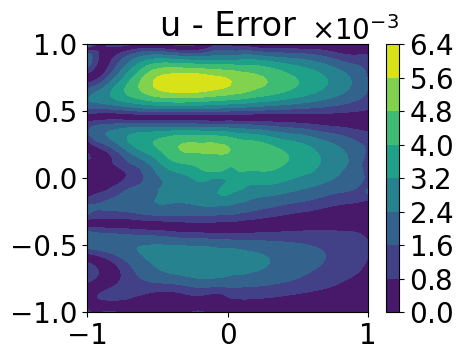}
        \caption{RFM (u)}
    \end{subfigure}
    \hfill
    \begin{subfigure}[t]{0.3\linewidth}
        \centering
         \includegraphics[width=\linewidth]{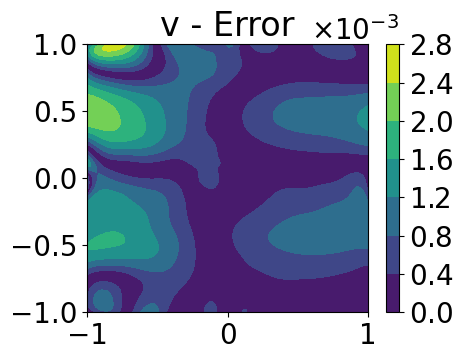}
        \caption{RFM (v)}
    \end{subfigure}
    \hfill
    \begin{subfigure}[t]{0.3\linewidth}
        \centering
        \includegraphics[width=\linewidth]{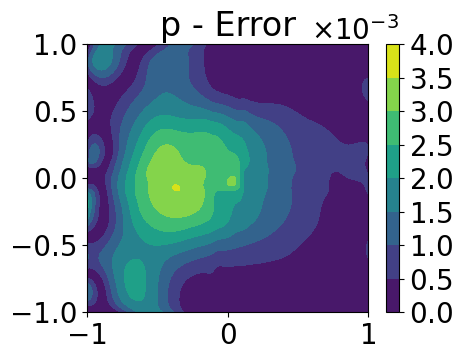}
        \caption{RFM (p)}
    \end{subfigure}

    \begin{subfigure}[t]{0.3\linewidth}
        \centering
        \includegraphics[width=\linewidth]{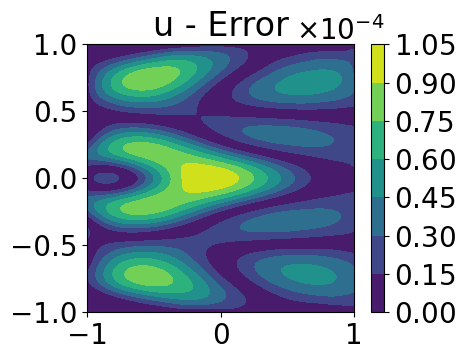}
        \caption{TransNet (u)}
    \end{subfigure}
    \hfill
    \begin{subfigure}[t]{0.3\linewidth}
        \centering
        \includegraphics[width=\linewidth]{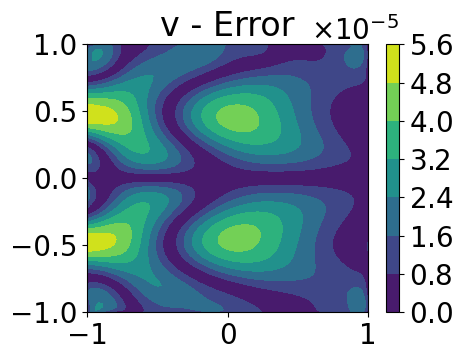}
        \caption{TransNet (v)}
    \end{subfigure}
    \hfill
    \begin{subfigure}[t]{0.3\linewidth}
        \centering
        \includegraphics[width=\linewidth]{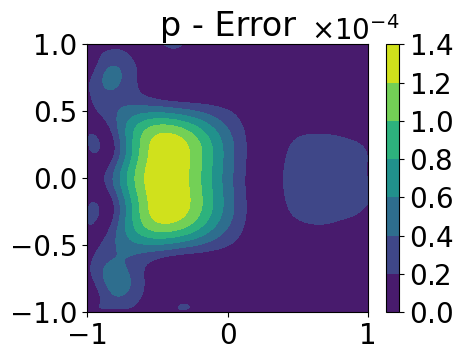}
        \caption{TransNet (p)}
    \end{subfigure}

    \caption{Pointwise absolute errors of $u,v,p$ across different methods (rows) for the steady-state Navier–Stokes equation.}
    \label{fig:error-comparison}
\end{figure}

In this experiment, $Re$ equals $40$, corresponding to $\nu = 0.025$. The computational domain is $\Omega = [-1, 1]^2$, and Dirichlet boundary conditions are prescribed by restricting the exact solution to the boundary $\partial\Omega$. To handle the nonlinearity of the Navier--Stokes equations, the Picard iteration is employed. At the $k$-th iteration, the convective term is linearized using the velocity field from the previous step, $\bm{u}_{\mathrm{NN}}^{k-1}$. The loss function minimized at iteration $k$ is:

\begin{equation}
    \mathcal{L}^{(k)} = \bm{u}_{\mathrm{NN}}^{k-1} \cdot \nabla \bm{u}_{\mathrm{NN}}^k + \nabla p_{\mathrm{NN}}^k - \nu \Delta \bm{u}_{\mathrm{NN}}^k.
\end{equation}
After 20 steps of Picard iteration, the pointwise absolute error of $u,v,p$ in different methods are shown in Figure~\ref{fig:error-comparison}.

\section{Discussion and Future Work}
\label{sec:discussion}

In this section, we provide a critical assessment of the proposed PD-OFM, focusing on its computational implications and the trade-offs between orthogonality and numerical performance.
\subsection{Efficiency}
\begin{itemize}
    \item \textbf{Regularizing Orthogonality}: The cost of computing the orthogonality loss is essentially negligible with respect to calculating PINN-loss, as it scales independently of the PDE's order and the neural network's depth.
    \item \textbf{Transferability}: Although the RFM benefits from zero-cost initialization, the pre-training cost of PD-OFM is a strategic investment: the learned nearly-orthogonal features exhibit strong transferability. This enables the reuse of discovered features, shifting the paradigm from random initializations to automated feature engineering.
\end{itemize}

\subsection{Accuracy}
\begin{itemize}
    \item \textbf{Gap in Theoretical Foundation of the Accuracy}: While imposing orthogonality consistently improves numerical accuracy in our experiments, a rigorous theoretical link remains an open question. Our results suggest that regularization expands the effective rank of the feature space (as indicated by the metrics in Section \ref{subsec:orthogonality}), but the theoretical mechanism by how this regularization reduces the approximation error needs further investigation.
    \item \textbf{Ill-conditioning}: Although PD-OFM decelerates the decay of singular values, the resulting system matrices may still exhibit some degree of ill-conditioning. Total elimination of feature correlation via gradient-based training remains difficult.
\end{itemize}

\section{Conclusion}

In this work, we introduced Physics-Driven Orthogonal Feature Method, a novel framework that integrates orthogonality regularization into PINNs to construct operator-adaptive orthogonal features. PD-OFM discovers features with greater representational capacity and improved transferability. Extensive numerical experiments on a range of PDEs including the Helmholtz, Poisson, Wave, and Navier-Stokes equations consistently demonstrate that PD-OFM significantly outperforms existing methods such as the random feature method \cite{Jingrun2022rfm} and TransNet \cite{ZHANG2024tarnsnet}. Overall, by leveraging orthogonal neural representations, PD-OFM bridges the gap between classical spectral methods and modern machine-learning-based solvers, offering a robust feature-learning framework for PDE applications.

\bmsection*{Data Availability Statement}

The data and source code that support the findings of this study are openly available in the GitHub repository at \url{https://github.com/kinjonfire/Physics-Driven-Discovery-of-Orthogonal-Features}. The code implements the PD-OFM framework and includes scripts to reproduce the numerical experiments presented in the paper.


\bmsection*{Funding}
D. Wang was partially supported by the National Natural Science Foundation of China (Grant No. 12422116),
Guangdong Basic and Applied Basic Research Foundation (Grant No. 2023A1515012199),
Shenzhen Science and Technology Innovation Program (Grant Nos. JCYJ20220530143803007, RCYX20221008092843046),
and Hetao Shenzhen-Hong Kong Science and Technology Innovation Cooperation Zone Project (No. HZQSWS-KCCYB-2024016).

\bibliography{references} 

\appendix

\section{Using Variational Loss to Solve Two-dimensional Poisson Equation on Different Domains}

This appendix studies the effect of using variational loss instead of traditional PINN loss to solve the two-dimensional Poisson equation on three types of domains: box, annulus, and L-shape. The equation is:
\begin{equation}
- \Delta u(x, y) = f(x, y), \quad (x, y) \in \Omega
\end{equation}
with zero Dirichlet boundary conditions. The variational loss \cite{yu2018deep} is used:
\begin{equation}
L_{\text{var}} = \frac{1}{2} \int_{\Omega} |\nabla u(x, y)|^2 \, dx \, dy - \int_{\Omega} f(x, y) u(x, y) \, dx \, dy.
\end{equation}

In practice, the variational loss is approximated by:
\begin{equation}
L_{\text{var}} \approx \frac{1}{2 N_q} \sum_{i=1}^{N_q} |\nabla u(x_i, y_i)|^2 - \frac{1}{N_q} \sum_{i=1}^{N_q} f(x_i, y_i) u(x_i, y_i)
\end{equation}
where $N_q$ is the number of quadrature points sampled in the domain. We use the neural network structure: two hidden layers with 500 neurons per layer and \texttt{tanh}$^3$ activation. The Adam optimizer is used for training. The proposed framework allows flexible design of PDE loss. The loss function is not limited to the classical PINN residual loss (strong form). Variational forms, weighted residual methods, weak forms, or other customized loss functions can also be applied within this structure. This flexibility opens up broad opportunities for future development in neural network-based PDE solvers.

\begin{figure}[h!]
    \centering
    \begin{subfigure}[t]{0.32\textwidth}
        \centering
        \includegraphics[width=\textwidth]{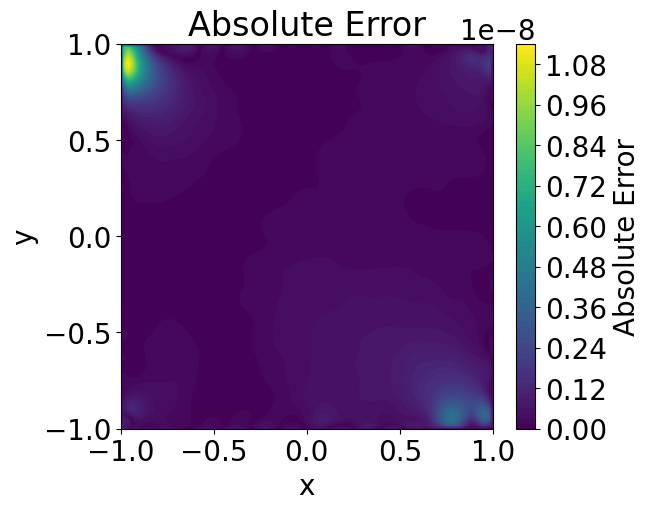}
        \caption{Box domain}
    \end{subfigure}
    \begin{subfigure}[t]{0.32\textwidth}
        \centering
        \includegraphics[width=\textwidth]{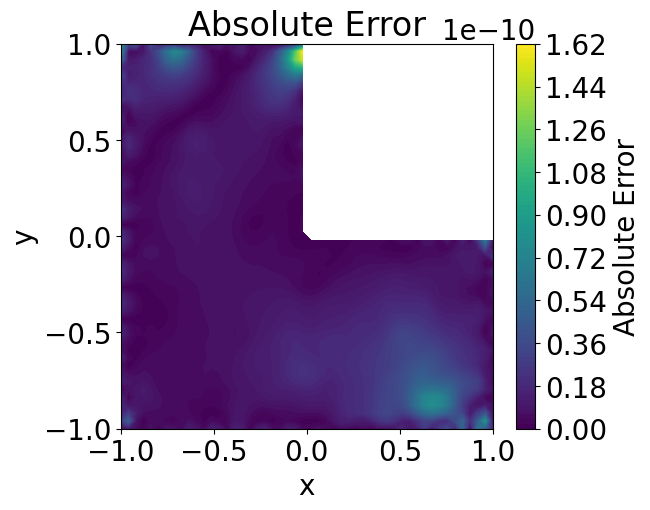}
        \caption{L-shape domain}
    \end{subfigure}
    \begin{subfigure}[t]{0.32\textwidth}
        \centering
        \includegraphics[width=\textwidth]{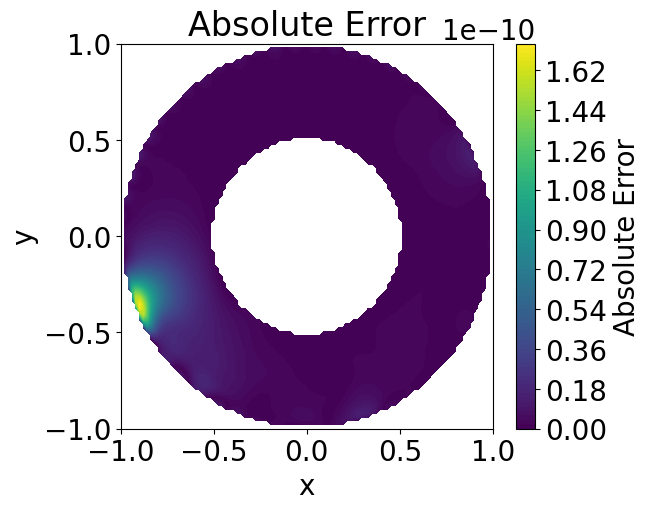}
        \caption{Annulus domain}
    \end{subfigure}
    \caption{Point-wise absolute error of variational loss solution on different domains.}
    \label{fig:error_domains}
\end{figure}

\section{Comparison of Transferability of PD-OFM and PD-FM}
\label{sec:comp}
\begin{table}[h]
\centering
\caption{Experiment 1: Error vs. Width and Depth in Box Domain Transferred Problem.}
\label{tab:exp1}
\begin{tabular}{@{}ccccccc@{}}
\toprule
\multirow{2}{*}{\textbf{Width}} &
\multicolumn{2}{c}{Depth = 1} &
\multicolumn{2}{c}{Depth = 2} &
\multicolumn{2}{c}{Depth = 3} \\
\cmidrule(lr){2-3}\cmidrule(lr){4-5}\cmidrule(lr){6-7}
& PD-OFM & PD-FM & PD-OFM & PD-FM & PD-OFM & PD-FM \\
\midrule
100 & $5.08\times10^{-2}$ & $5.13\times10^{-2}$ & $3.13\times10^{-2}$ & $5.66\times10^{-2}$ & $4.40\times10^{-2}$ & $6.87\times10^{-2}$ \\
200 & $1.08\times10^{-3}$ & $4.51\times10^{-3}$ & $4.94\times10^{-4}$ & $1.13\times10^{-3}$ & $1.13\times10^{-3}$ & $1.87\times10^{-3}$ \\
300 & $2.46\times10^{-4}$ & $2.49\times10^{-3}$ & $2.47\times10^{-5}$ & $3.06\times10^{-5}$ & $6.64\times10^{-5}$ & $4.96\times10^{-5}$ \\
400 & $8.13\times10^{-4}$ & $1.98\times10^{-3}$ & $6.34\times10^{-7}$ & $1.01\times10^{-6}$ & $7.66\times10^{-6}$ & $7.75\times10^{-7}$ \\
500 & $2.41\times10^{-4}$ & $6.28\times10^{-4}$ & $2.35\times10^{-8}$ & $2.12\times10^{-7}$ & $3.40\times10^{-7}$ & $7.66\times10^{-8}$ \\
600 & $7.85\times10^{-4}$ & $1.17\times10^{-3}$ & $2.73\times10^{-8}$ & $1.24\times10^{-7}$ & $3.16\times10^{-8}$ & $1.10\times10^{-8}$ \\
\bottomrule
\end{tabular}
\end{table}

\begin{table}[!h]
\centering
\caption{Experiment 2: Error vs. Width and Depth in L-shape Domain Transferred Problem.}
\label{tab:exp2}
\begin{tabular}{@{}ccccccc@{}}
\toprule
\multirow{2}{*}{\textbf{Width}} &
\multicolumn{2}{c}{Depth = 1} &
\multicolumn{2}{c}{Depth = 2} &
\multicolumn{2}{c}{Depth = 3} \\
\cmidrule(lr){2-3}\cmidrule(lr){4-5}\cmidrule(lr){6-7}
& PD-OFM & PD-FM & PD-OFM & PD-FM & PD-OFM & PD-FM \\
\midrule
100 & $8.64\times10^{-2}$ & $1.12\times10^{-1}$ & $6.76\times10^{-2}$ & $7.61\times10^{-2}$ & $8.53\times10^{-2}$ & $8.87\times10^{-2}$ \\
200 & $2.83\times10^{-3}$ & $1.09\times10^{-2}$ & $1.07\times10^{-3}$ & $1.98\times10^{-3}$ & $1.91\times10^{-3}$ & $2.58\times10^{-3}$ \\
300 & $1.34\times10^{-3}$ & $4.39\times10^{-3}$ & $4.86\times10^{-5}$ & $6.97\times10^{-5}$ & $1.19\times10^{-4}$ & $7.99\times10^{-5}$ \\
400 & $1.29\times10^{-3}$ & $4.33\times10^{-3}$ & $1.83\times10^{-6}$ & $3.74\times10^{-6}$ & $1.17\times10^{-5}$ & $1.95\times10^{-6}$ \\
500 & $1.13\times10^{-3}$ & $3.19\times10^{-3}$ & $8.32\times10^{-8}$ & $2.15\times10^{-6}$ & $5.87\times10^{-7}$ & $4.47\times10^{-7}$ \\
600 & $1.29\times10^{-3}$ & $2.52\times10^{-3}$ & $2.03\times10^{-8}$ & $9.30\times10^{-7}$ & $3.92\times10^{-8}$ & $1.39\times10^{-7}$ \\
\bottomrule
\end{tabular}
\end{table}

\begin{table}[!h]
\centering
\caption{Experiment 3: Error vs. Width and Depth in Annulus Domain Transferred Problem.}
\label{tab:exp3}
\begin{tabular}{@{}ccccccc@{}}
\toprule
\multirow{2}{*}{\textbf{Width}} &
\multicolumn{2}{c}{Depth = 1} &
\multicolumn{2}{c}{Depth = 2} &
\multicolumn{2}{c}{Depth = 3} \\
\cmidrule(lr){2-3}\cmidrule(lr){4-5}\cmidrule(lr){6-7}
& PD-OFM & PD-FM & PD-OFM & PD-FM & PD-OFM & PD-FM \\
\midrule
100 & $1.87\times10^{-2}$ & $2.41\times10^{-2}$ & $8.02\times10^{-3}$ & $1.29\times10^{-2}$ & $1.42\times10^{-2}$ & $2.11\times10^{-2}$ \\
200 & $2.12\times10^{-4}$ & $1.51\times10^{-3}$ & $4.96\times10^{-5}$ & $1.39\times10^{-4}$ & $1.38\times10^{-4}$ & $2.82\times10^{-4}$ \\
300 & $7.41\times10^{-5}$ & $6.34\times10^{-4}$ & $1.43\times10^{-6}$ & $2.04\times10^{-6}$ & $3.87\times10^{-6}$ & $3.65\times10^{-6}$ \\
400 & $7.28\times10^{-5}$ & $4.60\times10^{-4}$ & $1.99\times10^{-8}$ & $7.12\times10^{-8}$ & $2.49\times10^{-7}$ & $2.93\times10^{-8}$ \\
500 & $8.63\times10^{-5}$ & $3.80\times10^{-4}$ & $5.84\times10^{-10}$ & $9.02\times10^{-8}$ & $7.78\times10^{-9}$ & $2.41\times10^{-8}$ \\
600 & $8.08\times10^{-5}$ & $2.91\times10^{-4}$ & $5.62\times10^{-10}$ & $3.19\times10^{-8}$ & $3.56\times10^{-10}$ & $4.11\times10^{-9}$ \\
\bottomrule
\end{tabular}
\end{table}
\end{document}